\documentclass[12pt]{article}
\usepackage{amssymb}
\usepackage{enumerate}
\usepackage{amsmath}
\usepackage{graphicx} 
\usepackage{hyperref}
\usepackage{color}
\newtheorem{definition}{Definition}[section]

\newtheorem{algorithm}[definition]{Algorithm}

\DeclareMathAlphabet\mathbit
    \encodingdefault\rmdefault\bfdefault\itdefault
\DeclareOldFontCommand{\bi}{\normalfont\bfseries\itshape}{\mathbit}

\newcommand{\be}{\begin{equation}}
\newcommand{\ee}{\end{equation}}

\def\fakebold#1{\relax\ifvmode\leavevmode\fi%
\ifmmode%
\setbox0=\hbox{$#1$}%
\else%
\setbox0=\hbox{#1}%
\fi%
\kern-.02em\copy0 \kern-\wd0%
\kern .04em\copy0 \kern-\wd0%
\kern-.0125em\raise.02em\box0%
}%

\date{}
\begin{document}

\title{On the summation of divergent, truncated, and underspecified power series\\ via asymptotic approximants}
\author{Nathaniel S. Barlow$^1$, Christopher R. Stanton$^1$, Nicole Hill$^2$,\\ Steven J. Weinstein$^2$, and Allyssa G. Cio$^3$.\\
{\normalsize\itshape{$^1$School of Mathematical Sciences}},
{\normalsize\itshape{$^2$Department of Chemical Engineering}},\\
{\normalsize\itshape{$^3$Department of Industrial Engineering}}\\
 {\normalsize\itshape{Rochester Institute of Technology, Rochester, NY 14623, USA}}\\
}

\maketitle

This article has been accepted for publication in the Quarterly Journal of Mechanics and Applied Mathematics. The final published version can be found at https://academic.oup.com/qjmam/article/2895119/On-the-Summation-of-Divergent-Truncated-and

\begin{abstract} 
A compact and accurate solution method is provided for problems whose infinite power series solution diverges and/or whose series coefficients are only known up to a finite order.  The method only requires that either the power series solution or some truncation of the power series solution be available and that some asymptotic behavior of the solution is known away from the series' expansion point. Here, we formalize the method of \textit{asymptotic approximants} that has found recent success in its application to thermodynamic virial series where only a few to (at most) a dozen series coefficients are typically known. We demonstrate how asymptotic approximants may be constructed using simple recurrence relations, obtained through the use of a few known rules of series manipulation.   The result is an approximant that bridges two asymptotic regions of the unknown exact solution, while maintaining accuracy in-between.    A general algorithm is provided to construct such approximants.  To demonstrate the versatility of the method, approximants are constructed for three nonlinear problems relevant to mathematical physics: the Sakiadis boundary layer, the Blasius boundary layer, and the Flierl-Petviashvili monopole. The power series solution to each of these problems is underspecified since, in the absence of numerical simulation, one lower-order coefficient is not known;  consequently, higher-order coefficients that depend recursively on this coefficient are also unknown. 
The constructed approximants are capable of predicting this unknown coefficient as well as other important properties inherent to each problem. The approximants lead to new benchmark values for the Sakiadis boundary layer and agree with recent numerical values  for properties of the Blasius boundary layer and Flierl-Petviashvili monopole.  
\end{abstract}

\section{Introduction \label{sec:intro}}
Power series arise in virtually all applications of mathematical physics.  The  utility of such series is evident in the construction of approximations (e.g. finite differences) or as a rigorously determined solution to a problem.  Limitations generally inherent to power series solutions often inhibit their direct use, and they are more often useful in the implementation of numerical schemes.  For instance, a Taylor series representation of an unknown function may not converge, as it may have a finite radius of convergence arising from singularities (often complex) in the function it represents~\cite{Titchmarsh}.   Even when singularities are not a concern, higher-order terms of the series may be exceedingly difficult to compute~\cite{Masters}, which is especially problematic if the series converges slowly.

%\subsection{Asymptotic approximants}
Several techniques have been put forward to analytically continue and to accelerate convergence of divergent or slowly converging series; see, for example~\cite{Bender} (ch. 8),~\cite{BakerCritical} (chs. 14, 19, 20),~\cite{Chisholm},~\cite{Guttmann1},~\cite{Guttmann2},~\cite{BakerBook},~\cite{Andrianov},~\cite{FrostHarper}.  One of the more well-known techniques is the Pad\'e approximant method~\cite{BakerGammel}, which has similarities with the method of asymptotic approximants described herein.  A Pad\'e approximant (referred to as a Pad\'e) is the quotient of two polynomials, with $N$ total polynomial coefficients (distributed between the numerator and denominator), chosen such that the Taylor expansion of the Pad\'e reproduces the power series of interest up to $N$ terms.  This involves solving an algebraic system with $N$ unknowns (the Pad\'e coefficients).  After its construction, one may use the Pad\'e in lieu of the original series.  The $N^\text{th}$-order Pad\'e and $N^\text{th}$-order truncated original series will follow one another near the expansion point.  Ideally, the Pad\'e will then continue on beyond the radius of convergence (if one is present) and in general represent the actual solution better than the truncated series.  If chosen judiciously, a Pad\'e sequence (specified by a fixed difference in order between the numerator and denominator) may converge rapidly to the correct solution as $N$ is increased, whereas the original series may converge slowly or diverge.  In the case of a divergent series, the Pad\'e approximates the singularity that is presumed to be responsible for divergence~\cite{VanDyke}, and in doing so enables an accurate summation of the truncated series.

One drawback of the Pad\'e method is that it is not always clear beforehand which sequence is best suited for a given problem~\cite{VanDyke}.  If one has a power series alone, it is difficult to choose a correct Pad\'e form (see~\cite{Clisby,Guerrero,Tan} for example).  However, if the power series arises from a physical problem, it is likely that some additional conditions or context may be gleaned.  These conditions may be conveniently available, or may need to be independently derived.  For instance, if one has a power series representation of $f(x)$ and it is known that $f$ asymptotically approaches a constant value as $x\to\infty$, a ``symmetric'' Pad\'e sequence of equal order in the numerator and denominator will contain Pad\'es that all preserve this asymptotic condition, and this sequence will uniformly converge towards the correct solution as $N$ is increased.   
Baker and Gammel~\cite{BakerGammel} recognized this important result and went further to state:  If a series for $f(x)$ is known and $f\sim C x^p$ as $x\to\infty$ ($p$ being an integer), the uniformly convergent Pad\'e sequence is the one with a difference of $p$ between the numerator and denominator order, hence preserving the $x\to\infty$ behavior as $N$ is increased. 

Baker $\&$ Gammel's statement may be extended to asymptotic behaviors beyond integer power laws using approximants other than Pad\'es~\cite{BarlowJCP,BarlowAIChE,Barlow2015}.   An example of such an approximant is found in the review by Frost and Harper~\cite{FrostHarper}, used to find an approximate solution for the drag coefficient on a sphere in fluid-filled tube; the approximant they used incorporates the correct non-integer power-law asymptotic behavior.   The present work may be considered an extension of the perspective given in~\cite{FrostHarper}- namely the construction of approximants that are asymptotically consistent with known behavior in the vicinity of the domain boundaries while maintaining accuracy in-between.     In this paper, we aim to formalize this approach by defining \textit{asymptotic approximants} in general:
\begin{definition}
Given a power series representation of some function $f(x)$:
\begin{equation}
f=\sum_{n=0}^\infty a_n (x-x_0)^n,
\label{powerseries}
\end{equation}
and an asymptotic behavior 
\begin{equation}
f\sim C f_a(x)~\text{as}~x\to x_a,
\label{asymp}
\end{equation}
where $C$ is a constant, an \textit{asymptotic approximant} is any function $f_A(x)$ that may be expressed analytically in closed form and that satisfies the following three properties:
\begin{enumerate}
\item The $N$-term Taylor expansion of $f_A$ about $x_0$ is identical to the $N$-term truncation of~(\ref{powerseries}).
\item $\displaystyle{\lim_{x\to x_a} (f_A/f_a)=}$ constant for any $N$.
%\item If it is known that $f$ is continuous on the open interval between $x_0$ and $x_a$,  $f_A$ should also be continuous on this interval.  
\item The sequence of approximants converges for increasing $N$. 
\end{enumerate}
\label{defn}
\end{definition}
Choosing an approximant that satisfies the above definition will lead to a uniformly convergent sequence as $N$ is increased that preserves the correct asymptotic behavior.  Note that in the above definition, $f_a$ need not be exact; in fact, only the leading order is typically needed to construct an adequate approximant.  We seek approximants whose unknown coefficients can be generated with ease.  We do not wish for the construction of approximants to undergo more floating-point operations than performed in an accurate numerical solution.   Most of the approximant coefficients calculated here may be obtained recursively using a few simple series relations listed in Appendix~\ref{sec:formulae}.  This is in contrast with Pad\'es and their extensions, which often require the inversion of a linear system~\cite{Bender, BakerCritical,FrostHarper}.  As shall be apparent, a theme of this work is to find/impose the simplest asymptotic approximant form possible while maintaining the desired accuracy and precision.

Approximants are closed form functions that may be analytically differentiated to obtain auxiliary properties without the loss of accuracy that occurs with discretized solutions.  Also, physically relevant properties may be cast as unknowns within an approximant, and the approximant can be used as a predictor for such properties~\cite{Boyd1997,BarlowJCP,BarlowAIChE,Barlow2015}.  This feature of approximants provides a significant problem-solving advantage.

The paper is organized as follows.  In Section~\ref{sec:virial}, we review recent applications of asymptotic approximants towards the truncated virial series of thermodynamics.  The analysis is recast here as part of a new unified framework that only requires the use of a few simple series relations.  In Section~\ref{sec:algorithm}, an algorithm is provided to construct asymptotic approximants for problems in general.  The method is then applied to three nonlinear boundary value problems in Section~\ref{sec:bvp}. In Section~\ref{sec:Sakiadis}, the Sakiadis boundary layer problem is solved using two types of asymptotic approximants -  one which is both accurate and simple in form, and another which matches its higher-order asymptotic behavior and is capable of predicting the wall-shear coefficient and other properties to a higher precision than reported before. In Section~\ref{sec:Blasius},  the Blasius boundary layer problem is solved using a simple asymptotic approximant, which is also capable of predicting important quantities such as the wall-shear coefficient. In Section~\ref{sec:FP}, an asymptotic approximant is used to solve the Flierl-Petviashvili monopole problem and predict relevant properties, which are shown to agree with newly generated numerical results.   Key findings are summarized in Section~\ref{sec:summary}.

\section{Relevant background and results from virial-based approximants \label{sec:virial}}
%\subsubsection{Virial series}
In this section, we review the recent application of asymptotic approximants towards problems in thermodynamics.  In doing the review, simplifications and details not previously disclosed are elucidated here and also incorporated into the algorithm of Section~\ref{sec:algorithm}. 

The virial series is an equation of state formulated as a series expansion about the ideal gas limit~\cite{Mason} ($\rho\to 0$), and is expressed as  
\begin{equation}
P=kT\sum_{n=1}^NB_n(T)\rho^{n},\hspace{0.1in} B_1=1
\label{virial}
\end{equation}
where $P$ is the pressure, $\rho$ is the number density, $k$ is the Boltzmann constant, and $T$ is the temperature.  The virial coefficients $B_n$ are functions of $T$ as well as other physical parameters describing a given intermolecular potential.  Specifically, the $n^{th}$ virial coefficient is defined as an integral over the position of $n$ molecules~\cite{McQuarrie}.   In the absence of an exact equation of state, the virial series represents fluid properties at low density more accurately than any other theoretical or computational method, since it is effectively the Taylor series of the exact solution about $\rho$=0.  A barrier to its usage is that the number of integrals appearing for each viral coefficient increases rapidly and nonlinearly with the order of the coefficient.  Consequently, significant effort is spent to develop efficient algorithms to compute higher-order coefficients~\cite{Singh,Wheatley}.  As an indication of the difficulty involved, only twelve virial coefficients are currently known (to within some precision) for the hard-sphere model fluid~\cite{Wheatley}.  For more realistic model fluids, significantly fewer coefficients are known. 

Even with several coefficients, the virial series often converges only in a small region near $\rho=0$. Pad\'e approximants have been applied to extend this region, with a common example being the Carnahan-Starling equation of state for hard-sphere fluids~\cite{Carnahan}.  With only a finite number of terms available, the choice of approximant (be it Pad\'e or not), is important as convergence may be accelerated by choosing an approximant with consistent asymptotic behavior following Definition~\ref{defn}.  In the following subsections, we review the recent application of asymptotic approximants towards divergent and slowly converging truncated virial series to construct accurate isotherms and predict critical properties.  

\subsection{Soft-sphere fluids \label{sec:SS}}
\begin{figure*}[h!]
\centering
\includegraphics[width=2.6in]{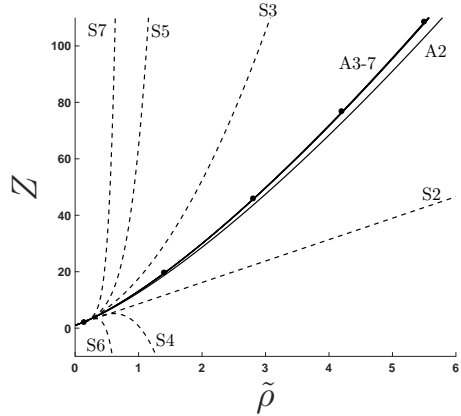}
\caption{Compressibility factor $Z$ versus reduced density $\tilde{\rho}$ for the $r^{-4}$ soft-sphere fluid.  A comparison is shown between the $N$-term virial series~(\ref{virial}) ($\---$) labeled as S$N$, the corresponding approximant~(\ref{SSAA}) ($\--$) labeled as A$N$, and data from molecular simulation~\cite{Rogers} ($\bullet$). The plot spans the entire fluid regime, as the last data point is where freezing occurs. Virial coefficients used to generate the curves are taken from~\cite{BarlowJCP}.}
\label{fig:softspheres}
\end{figure*}

To our knowledge, the first application of a non-Pad\'e asymptotic approximant towards virial series is given in~\cite{BarlowJCP} and is used to describe isotherms for soft-sphere fluids. This problem highlights the key steps in formulating an asymptotic approximant used later in this paper, and we thus review the technique carefully.   Soft-sphere fluids~\cite{Hansen} are defined by a molecular potential that behaves like $r^{-h}$, where $r$ is the distance between molecules and $h$ is a ``hardness'' parameter.  In the limit as $h\to\infty$, the soft-sphere molecular potential limits to that of the hard-sphere. For soft spheres, it is typical to rewrite~(\ref{virial}) in terms of reduced virial coefficients $\bar{B}_n$ and reduced density $\tilde{\rho}$ using the scalings provided in~\cite{BarlowJCP}.  Equation~(\ref{virial}) becomes
\begin{equation}
Z\equiv\frac{P}{\rho kT}=1+\sum_{n=2}^N\bar{B}_n(h)\tilde{\rho}^{~n-1},
\label{virialZ}
\end{equation}
where $Z$ is the compressibility factor and the coefficients $\bar{B}_n$ are now functions of the hardness $h$.  The above series diverges (see Fig.~\ref{fig:softspheres} for $h=4$), with a radius of convergence that decreases with decreasing $h$~\cite{BarlowJCP}. While~(\ref{virialZ}) accurately describes the low-density behavior, an approximant can be constructed by bridging~(\ref{virialZ}) with the large density limit for soft-spheres:
\begin{equation}
Z\sim C\tilde{\rho}^{~h/3} \text{ as } \tilde{\rho}\to\infty.   
\label{SS}
\end{equation}
Note that if $h/3$ is not an integer, it is impossible to capture the above behavior with a standard Pad\'e. However, an asymptotic approximant such as 
\begin{subequations}
\label{SSAA}
\begin{equation}
Z_A=\left(1+\sum_{n=1}^{N-1}A_n(h)\tilde{\rho}^{~n}\right)^{\frac{h/3}{N-1}}
\label{SSA}
\end{equation}
limits to~(\ref{SS}) for any $h$ and allows for an easy evaluation of its coefficients, such that the Taylor expansion of~(\ref{SSA}) about $\tilde{\rho}=0$ is exactly~(\ref{virialZ}).  An expression for the unknown coefficients $A_n$ is found by equating the original series~(\ref{virialZ}) with the approximant form~(\ref{SSA}) and then solving for the series in $A_n$:
 \[
\left(1+\sum_{n=2}^N\bar{B}_n(h)\tilde{\rho}^{~n-1}\right)^{\frac{N-1}{h/3}}=1+\sum_{n=1}^{N-1}A_n(h)\tilde{\rho}^{~n}.\]
It is now clear that the $A_n$ coefficients are the coefficients of the Taylor expansion of the left-hand side of the above expression about $\tilde{\rho}=0$.  Such coefficients are easily obtained using J. C. P. Miller's formula~\cite{Henrici} for recursively evaluating the expansion of a series raised to any power; see equation~(\ref{Miller}) in Appendix~\ref{sec:formulae}.   The recursion for the coefficients is
\begin{equation}
A_{n>0}=\frac{1}{n}\sum_{j=1}^{n}\left[\frac{j(N-1)}{h/3}-n+j\right]\bar{B}_{j+1}A_{n-j},~~A_0=1.
\label{SSC}
\end{equation}
\end{subequations}
Together,~(\ref{SSA}) and~(\ref{SSC}) provide an asymptotic approximant which converges for all $\tilde{\rho}$ in the fluid regime, as shown in~\cite{BarlowJCP} and Fig.~\ref{fig:softspheres} (for $h$=4).  Note that~(\ref{SSC}) reduces to the formulae for $A_2$ through $A_{10}$ given in the appendix of~\cite{BarlowJCP}, and can be used to generate higher-order approximants as higher-order virial coefficients become available.  The same procedure used above will be applied in Section~\ref{sec:bvp} to obtain asymptotic  approximants for  nonlinear boundary value problems of mathematical physics.  

\subsection{Prediction of Critical Properties \label{sec:critical}}
An important property of asymptotic approximants is their ability to predict unknown quantities.  Pad\'e and other approximants have long been used to predict critical properties of lattice models~\cite{Baker1961b,Thompson,Guttmann1,Guttmann2}.   With the advent of asymptotic approximants, this approach has recently yielded accurate predictions of critical properties of model fluids~\cite{BarlowAIChE,Barlow2015}.  The prediction begins by constructing the approximant in the same manner as the previous section.  

\begin{figure*}[h!]
\begin{tabular}{cc}(a)\hspace{-.05in}
\includegraphics[width=2.6in]{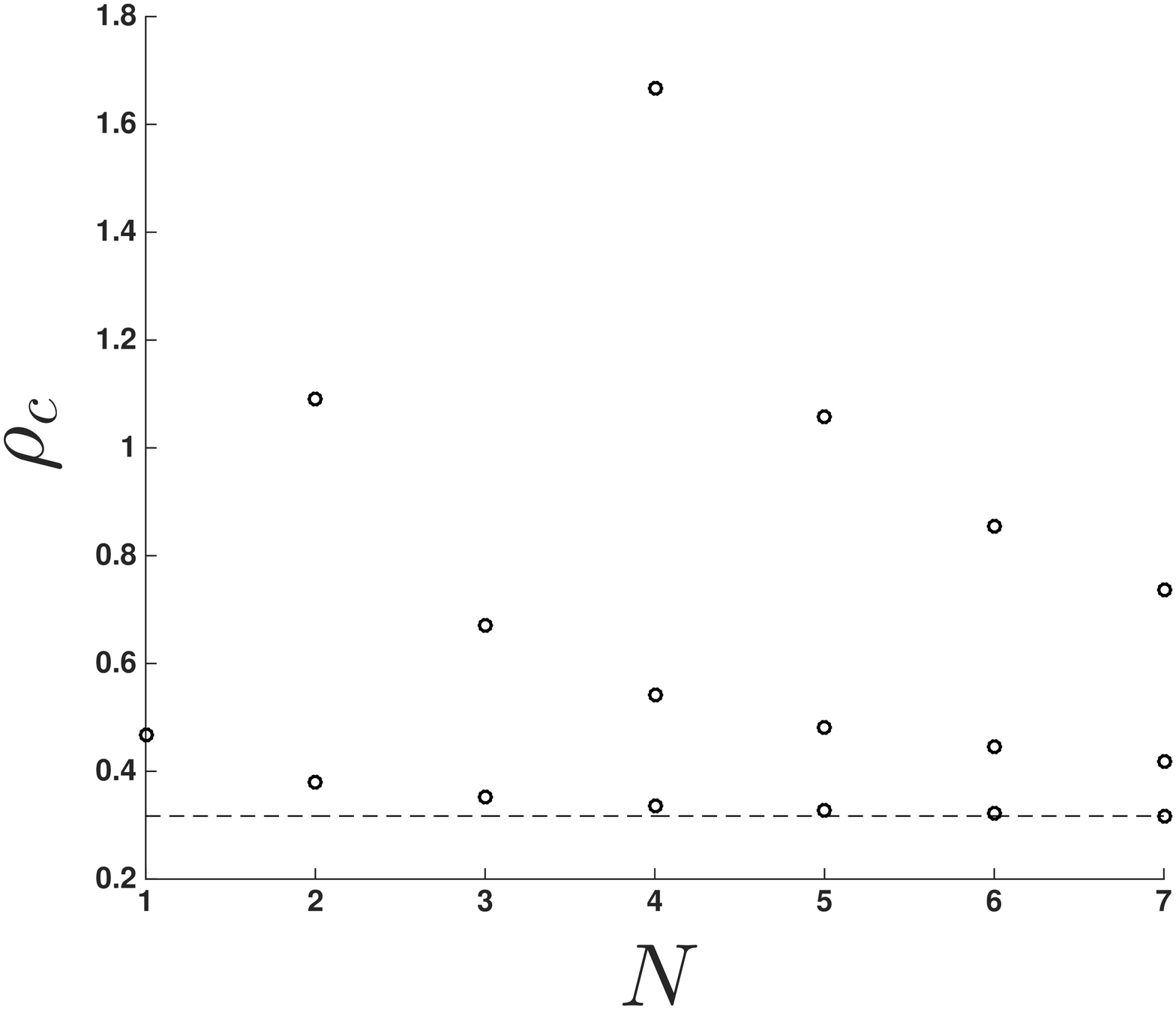}
(b)\hspace{-.05in}
\includegraphics[width=2.6in]{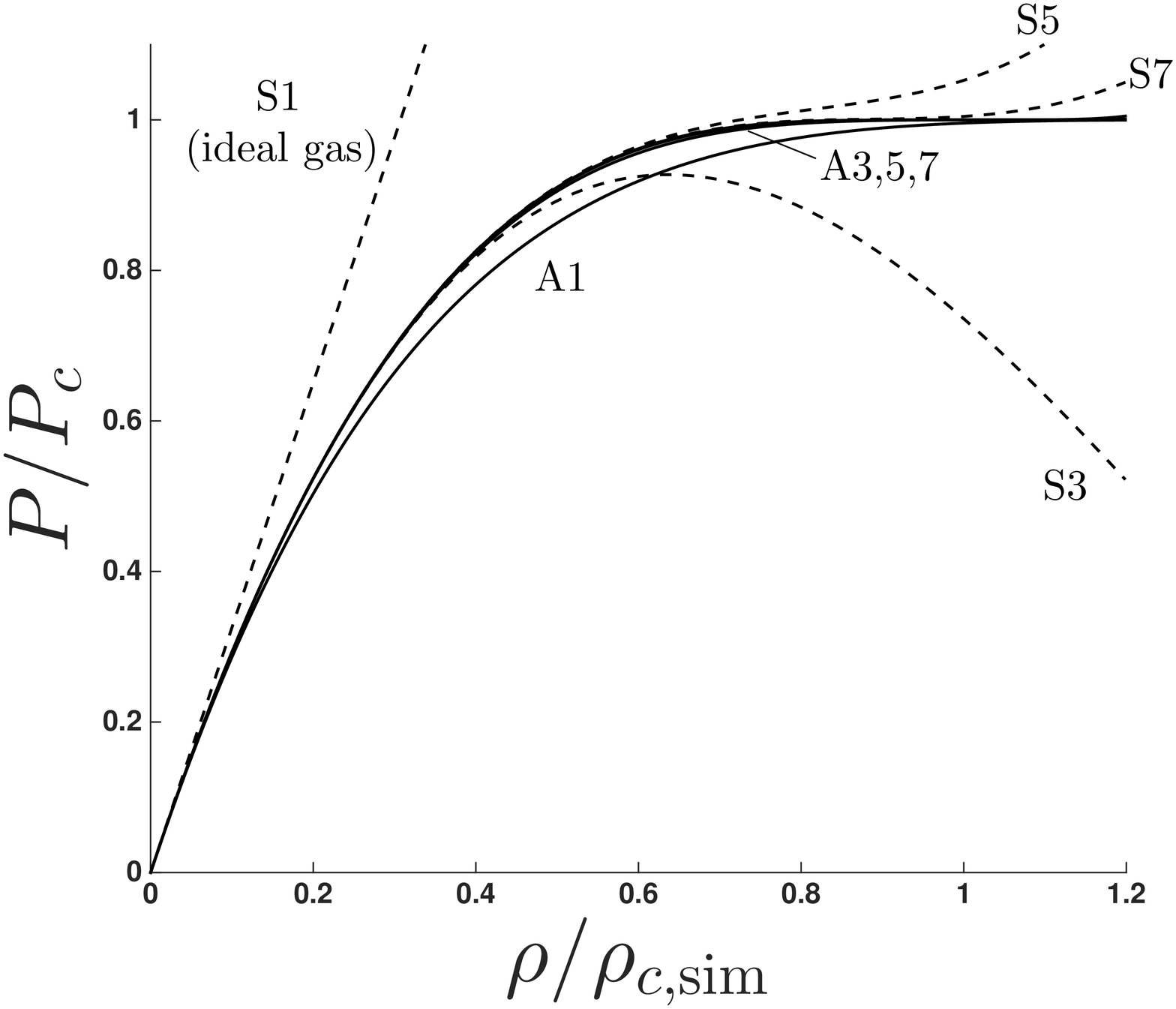}
\end{tabular}
\caption{(a) Predictions of the critical density from~(\ref{eq:rhoC}) for the Lennard-Jones fluid. A dashed line indicates the value predicted from Monte-Carlo simulations~\cite{Pellitero}.  (b)  Critical isotherm of the Lennard-Jones fluid. A comparison is shown between the $N$-term virial series~(\ref{virial}) ($\---$) labeled as S$N$ and the corresponding approximant~(\ref{C}) ($\--$) labeled as A$N$.  The most rapidly converging branch of $\rho_c$ roots (lowermost points) in (a) is used as an input to the approximants shown in (b).  Virial coefficients used to generate the plots are taken from~\cite{BarlowAIChE}.}
\label{fig:LJ}
\end{figure*}

The virial series describing $P(\rho)$ along the critical isotherm (i.e.~using coefficients evaluated at the critical temperature $T=T_c$) captures low density behavior but diverges due to a branch point singularity at the critical density $\rho_c$~\cite{Fisher}. The behavior of fluids near the critical point obeys a universal scaling law~\cite{Behnejad}
\begin{equation}
\left(\frac{P}{P_c}-1\right)\sim \pm D\left|\frac{\rho}{\rho_c}-1\right|^\delta,~T= T_c, ~\rho\rightarrow\rho_c,
\label{delta}
\end{equation}
where $P_c$ is the critical pressure, $D$ is a critical amplitude, and (for real fluids) $\delta$ is a positive non-integer critical exponent.  The knowledge of a low density series as well as an asymptotic behavior at higher density are the necessary ingredients for the application of an asymptotic approximant.  In particular, the series~(\ref{virial}) is accurate in a region where~(\ref{delta}) does not apply and inaccurate in the region where~(\ref{delta}) is accurate.  An approximant is provided in~\cite{BarlowAIChE} that captures both regions (low density and the critical region); it is given by
\begin{subequations}
\label{CC}
\begin{equation}
P_A=P_c-\sum_{n=0}^NA_n\rho^n\left(1-\frac{\rho}{\rho_c}\right)^\delta.
\label{CI}
\end{equation}
Following the same approach of the previous subsection, one may equate~(\ref{CI}) with the virial series~(\ref{virial}) evaluated at $T=T_c$ and solve for the series in $A_n$:
\[\left(P_c-kT_c\sum_{n=1}^NB_n(T_c)\rho^{n}\right)\left(1-\frac{\rho}{\rho_c}\right)^{-\delta}=\sum_{n=0}^NA_n\rho^n\]
where it becomes clear that the $A_n$ coefficients are, in fact, those in the Taylor expansion about $\rho=0$ of the left-hand side of the above expression.  The $A_n$ coefficients are obtained by replacing $(1-\rho/\rho_c)^{-\delta}$ with its expansion about $\rho=0$ and then taking the Cauchy product (equation~(\ref{Cauchy}) in Appendix A) of the two series in the left-hand side of the above expression.  This leads to Equation 4b of~\cite{BarlowAIChE}:
\begin{equation}
A_{n>0}=\frac{P_c \Gamma(\delta+n)}{n!\rho_c^n\Gamma(\delta)}-\frac{kT_c}{\Gamma(\delta)}\sum_{j=0}^{n-1}\frac{B_{n-j}\Gamma(\delta+j)}{\rho_c^j~j!},~A_0=P_c.
\label{CA}
\end{equation}
If $\rho_c$, $T_c$, and $P_c$ are known, one then inserts~(\ref{CA}) into~(\ref{CI}) to construct the approximant, which captures both the low-density and critical region, as shown in~\cite{BarlowAIChE}.  However, it is typically the case that critical properties are poorly known and difficult to compute through molecular simulation.  Fortunately, the approximant may be used to predict these quantities.  For instance, if $\rho_c$ is unknown and all other parameters are known, one may sacrifice the highest unknown coefficient $A_N$ and instead predict $\rho_c$.  This is equivalent to setting $A_N$=0 in~(\ref{CA}), leading to Equation 7 of~\cite{BarlowAIChE},
\begin{equation}
P_c-\frac{kT_cN!}{\Gamma(\delta+N)}\sum_{j=1}^N\frac{\Gamma(\delta+N-j)}{(N-j)!}B_j\rho_c^j=0
\label{eq:rhoC}
\end{equation}
\end{subequations}
which is a polynomial in $\rho_c$.  Although this predicts multiple roots for $\rho_c$, complex or negative values may be eliminated.   As one increases $N$ and tracks the roots, typically one set of positive roots will converge to a limit point more rapidly than others, as shown in Fig.~\ref{fig:LJ}a for the Lennard-Jones model fluid.  Taking the fastest converging $\rho_c$ sequence in the figure and substituting it into~(\ref{CA}) and~(\ref{CI}) for each $N$ leads to a set of approximant isotherms, shown in Fig.~\ref{fig:LJ}b.  Note in the figure that the approximant converges more rapidly and uniformly than the virial series as $N$ is increased.

Now if both $\rho_c$ and $P_c$ are unknown, one may sacrifice the two highest coefficients (i.e. $A_N=A_{N-1}=0$), leading to a system of 2 nonlinear equations to solve  for $\rho_c$ and $P_c$.  Recently, this approach has been used to predict properties and construct approximants that describe the entire critical region~\cite{Barlow2015} - not just the critical isotherm.  Results in Section~\ref{sec:bvp} demonstrate that this method is also useful in the prediction of quantities important to nonlinear boundary value problems.

As mentioned in the introduction, the virtue of approximants such as~(\ref{SSAA}) and (\ref{CC}) is the ease with which accurate derivative and integral properties are generated.  For example, the more comprehensive critical approximant given in~\cite{Barlow2015} for $P(\rho,T)$ can be integrated with respect to $\rho$ to accurately calculate the Helmholtz free energy, which may then be differentiated with respect to $T$ to find the internal energy and specific heat.   Even if these quantities are extracted numerically after an approximant is constructed, they may be resolved on a grid of arbitrary resolution with minimal computational expense (only memory usage), since the approximant is an analytic expression. 

\section{General procedure for constructing an asymptotic approximant \label{sec:algorithm}}
With the aim of formalizing the method of asymptotic approximants introduced in Section~\ref{sec:virial}, we provide the following algorithm.  This procedure is subsequently applied to power series solutions of nonlinear boundary-value problems in Section~\ref{sec:bvp}.
\begin{algorithm}
Steps for constructing an asymptotic approximant for a problem whose analytic solution $f(x)$ is elusive:
\begin{enumerate}[i.]
\item Find a power series solution of the given problem and truncate to $N$ terms: \[f_S=\sum_{n=0}^N a_n (x-x_0)^n.\]
Define $m$ as the number of unknowns embedded in the coefficients $a_n$; this value is used in step v below. Note that step i is independent of step ii and thus their order may be interchanged. 
\item Seek the asymptotic behavior $f\sim C f_a$ as $x\to x_a$, where $x_a\neq x_0$ and $C$ is a content.  This may come from either an asymptotic expansion about $x_a$, a boundary condition, or some independently known limiting behavior of $f$. Define $p$ as the number of unknowns associated with the asymptotic behavior; this value is used in step v below.
\item  Create an approximant function $f_A$ with $N+1$ unknowns $A_0\dots A_N$.  The function must (a) match the correct asymptotic behavior $f\sim C f_a$ as $x\to x_a$ to at least zeroth order and (b) be generalized to handle arbitrary $N$;  a function of an ``approximant series'' is well suited for this:
\[ f_A=f_A\left(\sum_{n=0}^{N} A_n [g(x)]^n\right).\]
This is the creative step.  The goal is to capture known features of $f$ while leaving flexibility to continuously merge two (potentially disparate) regions.  
\item Solve for $A_0\dots A_N$ by imposing the condition that the $N^{th}$-order Taylor expansion of $f_A$ about $x_0$ is equal to $f_S$.  If $g(x)=x$, this can be accomplished by setting $f_S=f_A$, isolating the approximant series within $f_A$ to one side and then Taylor expanding the other side using the relations in Appendix~\ref{sec:formulae}.  If $g(x)\neq x$, a matrix inversion may be required.  In either case, the resultant expressions for $A_0\dots A_N$ will be functions of the $m+p$ unknowns from steps i and ii.  
\item Solve the algebraic system $A_N=0,~A_{N-1}=0,~\dots~A_{N-m-p+1}=0$ for the $m+p$ unknowns from steps i and ii.  Choose the most rapidly converging and physically plausible roots (for example, the physical quantity the root represents may be real and positive).  In this step, unknown coefficients in the approximant series are sacrificed in order to solve for other unknowns embedded in these coefficients, while preserving the correct number of degrees of freedom. 
\item Evaluate all $A_n$ coefficients at the values obtained in step v and substitute these into the approximant $f_A$ from step iii. 
\item Ensure that $f_A$ converges as $N$ is increased.  As an additional metric of the approximant's ability to approach the correct asymptotic behavior, one may plot $f_A/f_a$ and confirm that, for increasing $N$, this ratio converges to a constant as $x\to x_a$. 
\item Consistency Checks:  Expand $f_A$ about $x_0$ to $N^{th}$ order; if the steps were executed correctly this will match the first $N$ terms of $f_S$.  Take the limit of $f_A$ as $x\to x_a$ and confirm that it agrees with $f_a$ to zeroth order.   Finally, if available, verify that the approximant agrees with a full numerical solution of the problem.

% where $m$ represents the number of unknowns are already embedded in the $a_n$ coefficients
   
\end{enumerate}
\end{algorithm}

%The problems are presented in an order based on the quality of the approximant used, from increasing to decreasing. 

\section{Application to nonlinear boundary value problems \label{sec:bvp}}
Asympotic approximants are now applied to three nonlinear boundary value problems. In Section~\ref{sec:Sakiadis}, we find approximant solutions to the Sakiadis boundary layer problem, where we are able to take advantage of the full asymptotic behavior, enabling us to obtain new benchmark values.  In Section~\ref{sec:Blasius}, an approximant is constructed for the Blasius boundary layer problem using the leading order asymptotic behavior; this approximant confirms existing numerically-obtained benchmark values. In Section~\ref{sec:FP}, an approximant for the Flierl-Petviashvili monopole problem is constructed using the leading order asymptotic behavior.  In this application, newly obtained numerical results are confirmed and a new asymptotic constant is explored. 

\subsection{The Sakiadis problem \label{sec:Sakiadis}}
The Sakiadis problem describes steady-state developing flow created by a moving plate in an otherwise stagnant fluid~\cite{Sakiadis}.   It shares the same boundary layer governing equation as that of the Blasius problem (discussed in Section~\ref{sec:Blasius}), but with different boundary conditions.  As shown in~\cite{Sakiadis}, a similarity transform of these equations leads to the following nonlinear boundary value problem in $f(\eta)$:
\begin{subequations}
\begin{equation}
2f'''+ff''=0
\label{sDE}
\end{equation}
\begin{equation}
f(0)=0, ~ f'(0)=1, ~ f'(\infty) = 0.
\label{sBC}
\end{equation}
\label{eq:Sakiadis}
\end{subequations}
As shown by Blasius~\cite{Blasius}, substituting the solution form
\begin{subequations}
\begin{equation}
f= \sum\limits_{n=0}^\infty a_n\eta^n
\label{power series}
\end{equation}
into~(\ref{sDE}) leads to the following recursion for the series coefficients
\begin{equation}
a_{n+3} =  \frac{-\sum\limits_{j=0}^n (j+1)(j+2)~a_{j+2}~a_{n-j}}{2(n+1)(n+2)(n+3)}.
\label{coefficients}
\end{equation}
\label{eq:series}
\end{subequations}
The recursion above  requires knowledge of the first three coefficients in order to generate the full series.   The first two coefficients $a_0=f(0)$ and $a_1=f'(0)$ are given by the first two boundary conditions in~(\ref{sBC}). The third coefficient is prescribed by the ``wall shear'' parameter \[f''(0)\equiv\kappa,\] and is typically obtained from a numerical solution of (\ref{eq:Sakiadis}) (see ~\cite{Cortell,Eftekhari,Fazio} for example). Recently, a semi-analytical fixed point method was used to estimate $\kappa$; see~\cite{Xu} for details and a comprehensive review and comparison of previous numerical predictions.   Here, asymptotic approximants will be used to predict $\kappa$. 

From the infinite condition in~(\ref{sBC}), it follows that the solution limits to a constant $C$ at large $\eta$:
\begin{equation}
\lim_{\eta\to\infty}f\equiv C.
\label{C}
\end{equation}
Higher order asymptotic methods provided in Appendix~\ref{sec:SakiadisA} enable the identification of an additional asymptotic constant, $G$, as:
\begin{equation}
\lim_{\eta\to\infty}\left[(f-C)e^{\eta C/2}\right]\equiv G.
\label{G}
\end{equation}
The constants $G$, $C$, and $\kappa$ will be predicted using asymptotic approximants in the following subsections.

\subsubsection{Simple asymptotic approximant \label{sec:Sakiadis1}}
First, we consider a simple approximant of the form 
 \begin{subequations}
 \label{SakiadisA}
 \begin{equation}
% f_A=C-\frac{C}{\displaystyle1+\sum_{n=1}^NA_n~x^n},
  f_A=C-C\left(1+\sum_{n=1}^{N}A_n~\eta^n\right)^{-1},
 \label{SakiadisApproximant}
 \end{equation}
 which automatically satisfies both~(\ref{sBC}) and~(\ref{C}).  The assumed form~(\ref{SakiadisApproximant}) has the appearance of a Pad\'e, but in fact it is not. Since $C$ is an unknown to be predicted by the approximant, the above cannot be considered a Pad\'e for the quantity $f-C$.  Also, if the terms in~(\ref{SakiadisApproximant}) were combined through a common denominator, a specific relationship between the numerator and denominator coefficients would be obtained.  Thus~(\ref{SakiadisApproximant}) cannot be considered a Pad\'e for $f$.  We shall see below that the approximant given by~(\ref{SakiadisApproximant}) contains coefficients which are, in fact, easier to compute than those of a standard Pad\'e. 
 
 The coefficients $A_0\dots A_{N}$,  $C$, and $\kappa$ are calculated such that the $N$-term Taylor expansion of~(\ref{SakiadisApproximant}) about $\eta=0$ is exactly equal to the $N$-term truncation of~(\ref{eq:series}).  Equating~(\ref{SakiadisApproximant}) with~(\ref{eq:series}) and re-arranging leads to the following
\[\left[1-\frac{1}{C}\sum\limits_{n=0}^N a_n\eta^n\right]^{-1}=1+\sum_{n=1}^{N}A_n~\eta^n.\]
Expanding the left-hand side above (utilizing equation~(\ref{inverse}) in Appendix~\ref{sec:formulae}) provides the following recursion for the coefficients, which are now functions of $\kappa$ and $C$:
\begin{equation}
A_{n>0}=\frac{1}{C}\sum_{j=1}^na_{j}~A_{n-j},~~A_0=1.
\label{SakiadisCoefficients}
\end{equation}
We now sacrifice the coefficients $A_N$ and $A_{N-1}$ to simultaneously predict $\kappa$ and $C$.  Setting  $A_N=A_{N-1}=0$ in~(\ref{SakiadisCoefficients}), we arrive at 
\begin{align}
%\nonumber
\sum_{j=1}^Na_jA_{N-j}=0,~~~
\sum_{j=1}^{N-1}a_jA_{N-1-j}=0.
\label{kappaeqn}
\end{align}
\end{subequations}
Once the preceding coefficients are obtained from~(\ref{SakiadisCoefficients}) and substituted into the above,~(\ref{kappaeqn}) becomes a system of nonlinear equations in $\kappa$ and $C$, and may be solved algebraically.  The most rapidly converging $\kappa$ and $C$ roots of~(\ref{kappaeqn}) are given in Table~\ref{table:Skappa} up to an optimal asymptotic truncation of  $N=11$.  That is, for $N>11$, $|\kappa_{N+1}-\kappa_N|$ no longer consistently decreases as $N$ is increased.   

We also record the value $S$ as the magnitude of the singularity, $\eta_s$, closest to $\eta$=0 in the complex $\eta$ plane for the Sakiadis problem, as predicted by~(\ref{SakiadisA}) and listed in Table~(\ref{table:Skappa}).  This value limits the radius of convergence of the series~(\ref{eq:series}).

%The coefficients $A_0\dots A_{N-2}$  $C\equiv f(\infty)$, and $\kappa$ are calculated such that the $N$-term Taylor expansion of~(\ref{SakiadisApproximant}) about $\eta=0$ is exactly equal to~(\ref{eq:series}).
{\renewcommand{\arraystretch}{1.1}
\begin{table}[h!]
\centering
\caption{Predictions from approximants for the Sakiadis problem.  For comparison, the numerical values of $\kappa$ and $C$ in the bottom row were computed using the Chebfun~\cite{chebfun} package to solve~(\ref{eq:Sakiadis}) via rectangular collocation with $\eta=\infty$ replaced with a finite surrogate $\eta_\infty$ given in the table. Note that the recent $\kappa$ estimates given in~\cite{Andersson,Xu} (reported to 7 digits) agree with the first 7 digits of our numerical predictions and predictions from~(\ref{Exp}).    The numerical values of $S$ were obtained by inserting the numerical $\kappa$ and $C$ values into~(\ref{SakiadisA}), locating the singularity of smallest magnitude, and systematically increasing the approximant order (up to $N=200$), where the value of $S$ has converged to within the reported digits.}
%\small
\label{table:Skappa}
\setlength{\tabcolsep}{.3em}
\footnotesize
\begin{tabular}{ccccccc}
\hline
$N$       & \begin{tabular}{c}$\kappa$\\from~(\ref{SakiadisA})\end{tabular} & \begin{tabular}{c}$C$\\ from~(\ref{SakiadisA})\end{tabular}   & \begin{tabular}{c}$S$\\ from~(\ref{SakiadisA})\end{tabular} & \begin{tabular}{c}$\kappa$\\ from~(\ref{Exp})\end{tabular}   & \begin{tabular}{c}$C$\\ from~(\ref{Exp})\end{tabular}      & \begin{tabular}{c}$G$\\ from~(\ref{Exp})\end{tabular}     \\ \hline\hline
5         & -0.3879       & 2.1607     & 3.76     & -0.464241808775 & 1.57791591603 & -1.9379232426 \\
7         & -0.4341       & 1.7275     & 3.98     & -0.445957694266 & 1.61001774844 & -2.0928590899 \\
9         & -0.4421       & 1.6418     & 4.05     & -0.443970941603 & 1.61533742596 & -2.1255327825 \\
11        & -0.4430       & 1.6284     & 4.07     & -0.443769095918 & 1.61603802403 & -2.1306162759 \\
13        & $\--$         & $\--$      & $\--$    & -0.443750163304 & 1.61611659263 & -2.1312640887 \\
15        & $\--$         & $\--$      & $\--$    & -0.443748473247 & 1.61612459984 & -2.1313373822 \\
17        & $\--$         & $\--$      & $\--$    & -0.443748326909 & 1.61612536888 & -2.1313450758 \\
19        & $\--$         & $\--$      & $\--$    & -0.443748314499 & 1.61612543984 & -2.1313458428 \\
21        & $\--$         & $\--$      & $\--$    & -0.443748313462 & 1.61612544619 & -2.1313459165 \\
23        & $\--$         & $\--$      & $\--$    & -0.443748313376 & 1.61612544675 & -2.1313459234 \\
25        & $\--$         & $\--$      & $\--$    & -0.443748313368 & 1.61612544681 & -2.1313459241 \\
26        & $\--$         & $\--$      & $\--$    & -0.443748313384 & 1.61612544669 & -2.1313459226 \\
27        & $\--$         & $\--$      & $\--$    & -0.443748313370 & 1.61612544680 & -2.1313459239 \\
28        & $\--$         & $\--$      & $\--$    & -0.443748313369 & 1.61612544680 & -2.1313459240 \\
29        & $\--$         & $\--$      & $\--$    & -0.443748313369 & 1.61612544681 & -2.1313459241 \\
30        & $\--$         & $\--$      & $\--$    & -0.443748313369 & 1.61612544681 & -2.1313459241 \\  \hline 
numerical & $\kappa$ & $C$ & \begin{tabular}{c}$S$\\ from~(\ref{SakiadisA})\end{tabular} & $\eta_\infty$ &&\\
\hline\hline
shooting~\cite{Cortell} & -0.44374733 & 1.61611200 & 4.07217&20& & \\
quadrature~\cite{Eftekhari}& -0.44374831 & 1.61612518 & 4.07217&23.20512 &  &  \\ 
collocation& -0.44374831   & 1.61612545 & 4.07217  & 30     &     &  \\ \hline
\end{tabular}\end{table}

Lastly, $C$ and $\kappa$ are substituted into~(\ref{SakiadisCoefficients}) for all non-zero coefficients ($A_1$ through $A_{N-2}$) to construct the approximant series in~(\ref{SakiadisApproximant}). The approximant~(\ref{SakiadisA}) is compared with the series and numerical solution of~(\ref{eq:Sakiadis}) in Fig.~\ref{fig:sakiadisCrude}. Although we are using a simple approximant form that enables correspondingly straightforward coefficient generation, the approximant converges surprisingly well, as can be seen in Fig.~\ref{fig:sakiadisCrude}a.  In Fig.~\ref{fig:sakiadisCrude}b, we confirm that $f'$ and $f''$ are accurately obtained from the approximant, once $f$ has reasonably converged (here, for $N=11$).  Unfortunately, our estimates for $\kappa$ and $C$ cease to converge for $N>11$, which limits our ability to predict the asymptotic constant $G$ given by~(\ref{G}).  This can been seen in a plot of the effective $G$
\begin{equation}
G_\text{eff}=(f-C)e^{\eta C/2},~G_\text{eff}\to G~\text{as}~\eta\to\infty
\label{Geff}
\end{equation}
versus $\eta$, which, when computed by the approximant~(\ref{SakiadisA}), only converges in a small region near $\eta=0$ as shown in Fig.~\ref{fig:sakiadisCrude}b.  The reason for this poor convergence behavior is that the approximant form~(\ref{SakiadisA}) does not incorporate the asymptotic correction behavior that leads to~(\ref{G}), and thus $G_\text{eff}$ cannot converge uniformly. 

\begin{figure*}[h!]
\begin{tabular}{cc}(a)\hspace{-.05in}
\includegraphics[width=2.6in]{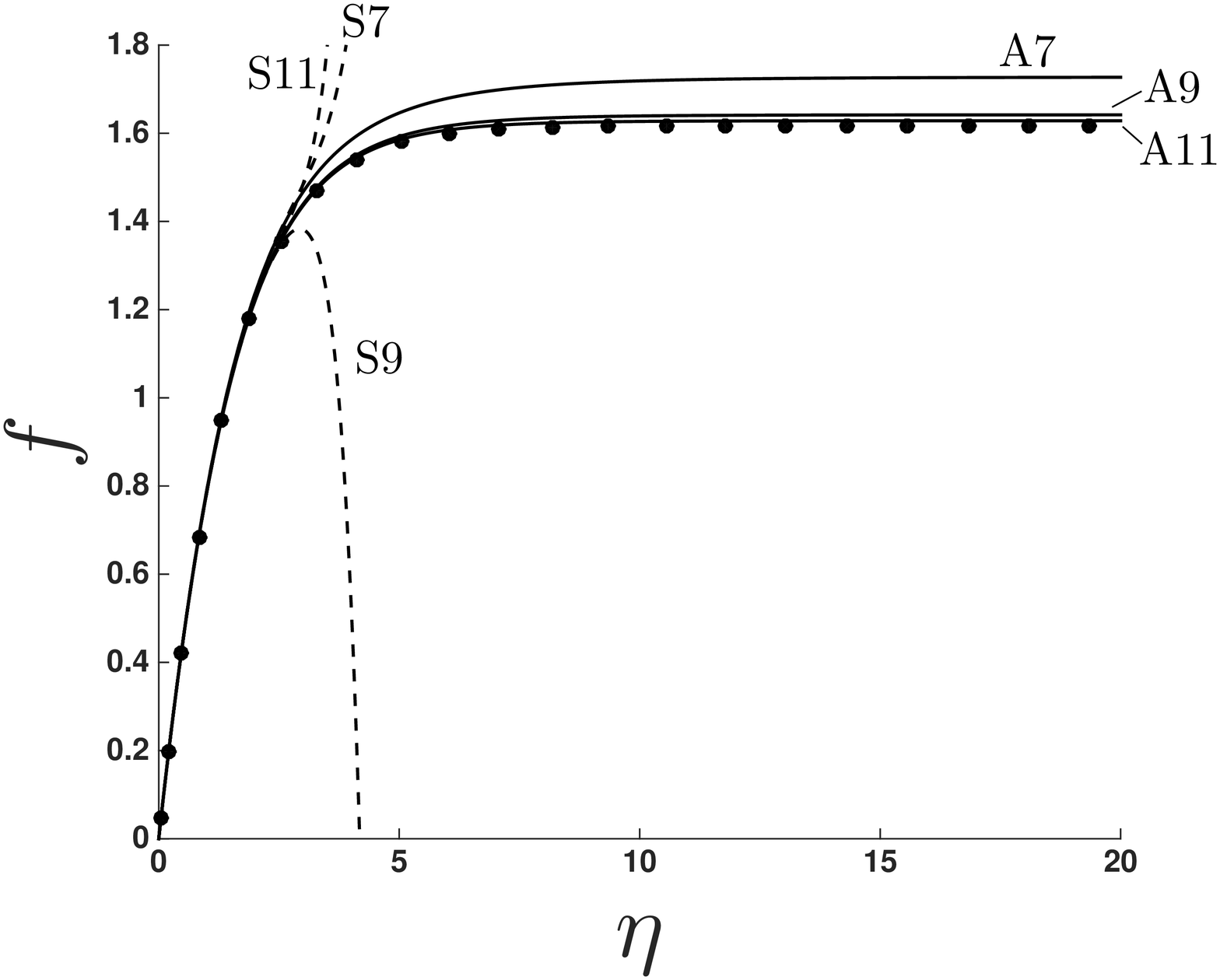}
(b)\hspace{-.05in}
\includegraphics[width=2.6in]{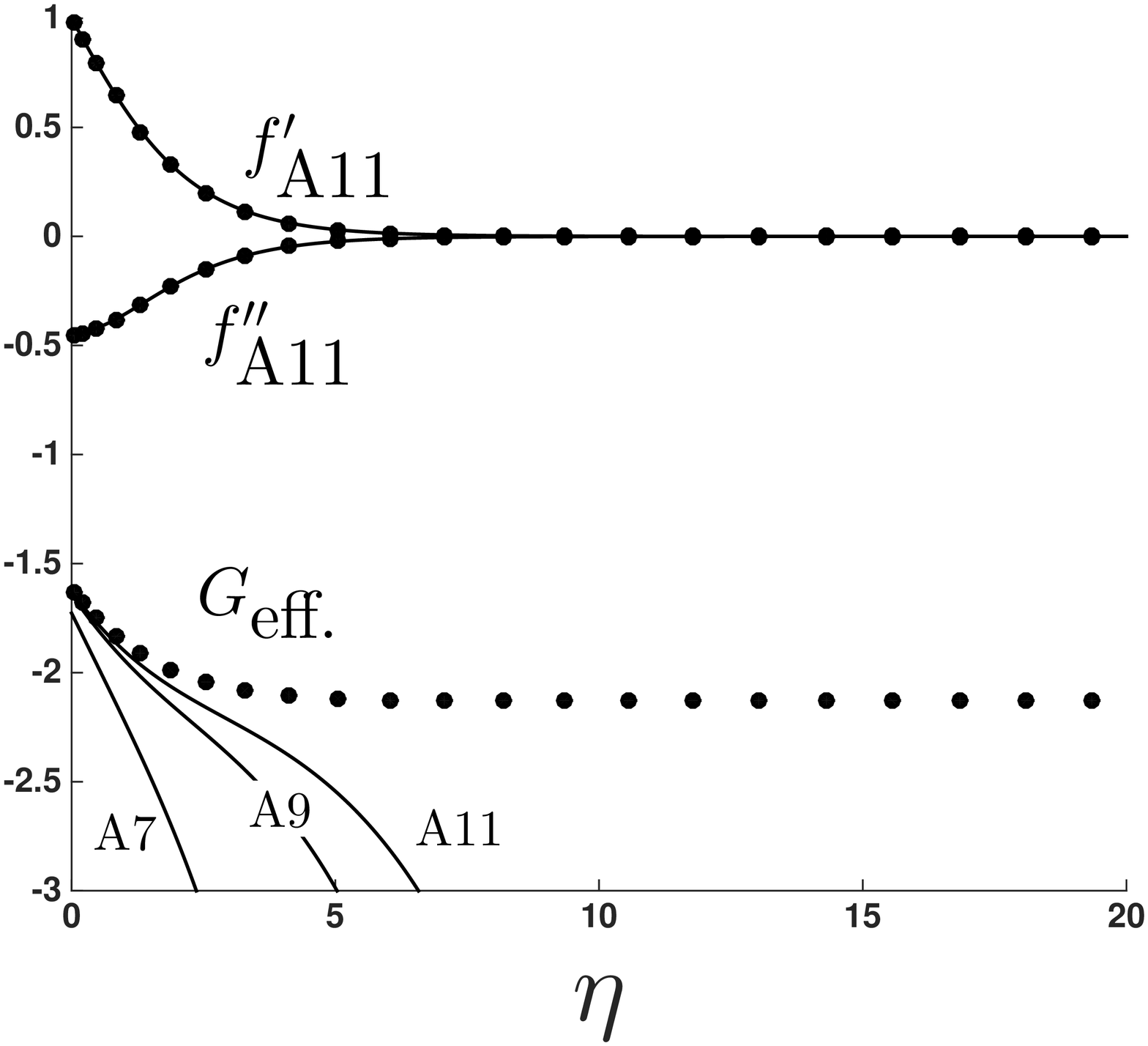}
\end{tabular}
\caption{(a) Comparison between the $N$-term series solution~(\ref{eq:series}) ($\---$) labeled as S$N$, the corresponding approximant~(\ref{SakiadisA}) ($\--$) labeled as A$N$, and the numerical solution~\cite{Cortell} ($\bullet$) to the Sakiadis problem~(\ref{eq:Sakiadis}). (b) Derivatives of approximant~(\ref{SakiadisA}) with $N$=11 and $G_\text{eff}$ given by~(\ref{Geff}) compared with the numerical solution~\cite{Cortell} ($\bullet$).}
\label{fig:sakiadisCrude}
\end{figure*}

\subsubsection{Asymptotic approximant with exponential corrections \label{sec:Sakiadis2}}
We now consider an approximant, 
\begin{subequations}
\label{Exp}
\begin{equation}
f_A=C+\sum_{n=1}^{N}A_ne^{-nC\eta/2},
\label{SakiadisExp}
\end{equation}
that not only matches the zeroth order infinite behavior, but also the form of the exponential corrections from the asymptotic expansion of $f$ as $\eta\to\infty$ provided in Appendix~\ref{sec:SakiadisA}.  Note that $A_1$ in the series above is the unknown asymptotic constant $G$ defined in~(\ref{G}).   After expanding~(\ref{SakiadisExp}) about $\eta=0$ and matching like terms with those of the series solution ~(\ref{eq:series}), we arrive at a system of $N$ linear equations,

\begin{equation}
\left[ \begin{array}{ccccc}

1^0 & 2^0 & 3^0 & \cdots & N^0 \\ 1^1 & 2^1 & 3^1 & \cdots & N^1 \\ 1^2 & 2^2 & 3^2 & \cdots & N^2 \\ \vdots & \vdots & \vdots & \vdots & \vdots \\ 1^{N-1} & 2^{N-1} & 3^{N-1} & \cdots & N^{N-1}\end{array} \right]\left[ \begin{array}{c}

A_1\\ A_2 \\ A_3 \\ \vdots \\ A_N\end{array} \right]=\left[ \begin{array}{c}

0!~a_0-C\\ 1!~(-2/C)~a_1 \\ 2!~ (-2/C)^2 ~a_2 \\ \vdots \\ (N-1)!~ (-2/C)^{N-1} ~a_{N-1} \end{array} \right],
\label{matrix}
\end{equation}
where the coefficient matrix is a Vandermonde matrix with a known explicit formula for its inverse~\cite{Vandermonde}. Once the system~(\ref{matrix}) is solved, the $A_n$ coefficients become functions of the unknowns $C$ and $\kappa$.  Like the previous approximant, we sacrifice the last two coefficients and solve the nonlinear equations 
\begin{equation}
A_N(\kappa,C)=0,~A_{N-1}(\kappa,C)=0
\end{equation}
\end{subequations}
 to obtain values for $C$ and $\kappa$, provided in the 5$^\text{th}$ and 6$^\text{th}$ columns of Table~\ref{table:Skappa}.  In the table, we report values up to $N=30$, after which the 13$^\text{th}$ digit oscillates without converging to a discernible  limit point.  The converged values \[ \kappa=-0.443748313369\] and \[C=1.61612544681\] obtained from approximant~(\ref{Exp}) agree with numerical predictions, and provide 5 additional digits of precision.  Also, for the first time, a value of the asymptotic constant \[G=-2.1313459241\] is reported; convergence to this value is shown in the rightmost column of Table~\ref{table:Skappa}, obtained by substituting the predictions for $C$ and $\kappa$ into the expression for $A_1$.   Once $C$ and $\kappa$ are substituted into the expressions (arising from the solution to~(\ref{matrix})) for all non-zero coefficients ($A_1$ through $A_{N-2}$), the approximant series in~(\ref{SakiadisExp}) may be constructed. 

The approximant~(\ref{Exp}) is compared with the series and numerical solution of~(\ref{eq:Sakiadis}) in Fig.~\ref{fig:sakiadis}.   The new approximant converges more rapidly than the simpler approximant (see Fig.~\ref{fig:sakiadisCrude} for comparison).  Also, the new approximant captures the exponential corrections inasmuch as $G_\text{eff}$ converges into the large $\eta$ asymptotic regime and approaches $G$, as shown in Fig.~\ref{fig:sakiadis}b.   This behavior is expected, as the correct asymptotic behavior has been incorportated explicitly into the approximant and convergence should be uniform for large $
\eta$.

\begin{figure*}[h!]
\begin{tabular}{cc}(a)\hspace{-.05in}
\includegraphics[width=2.6in]{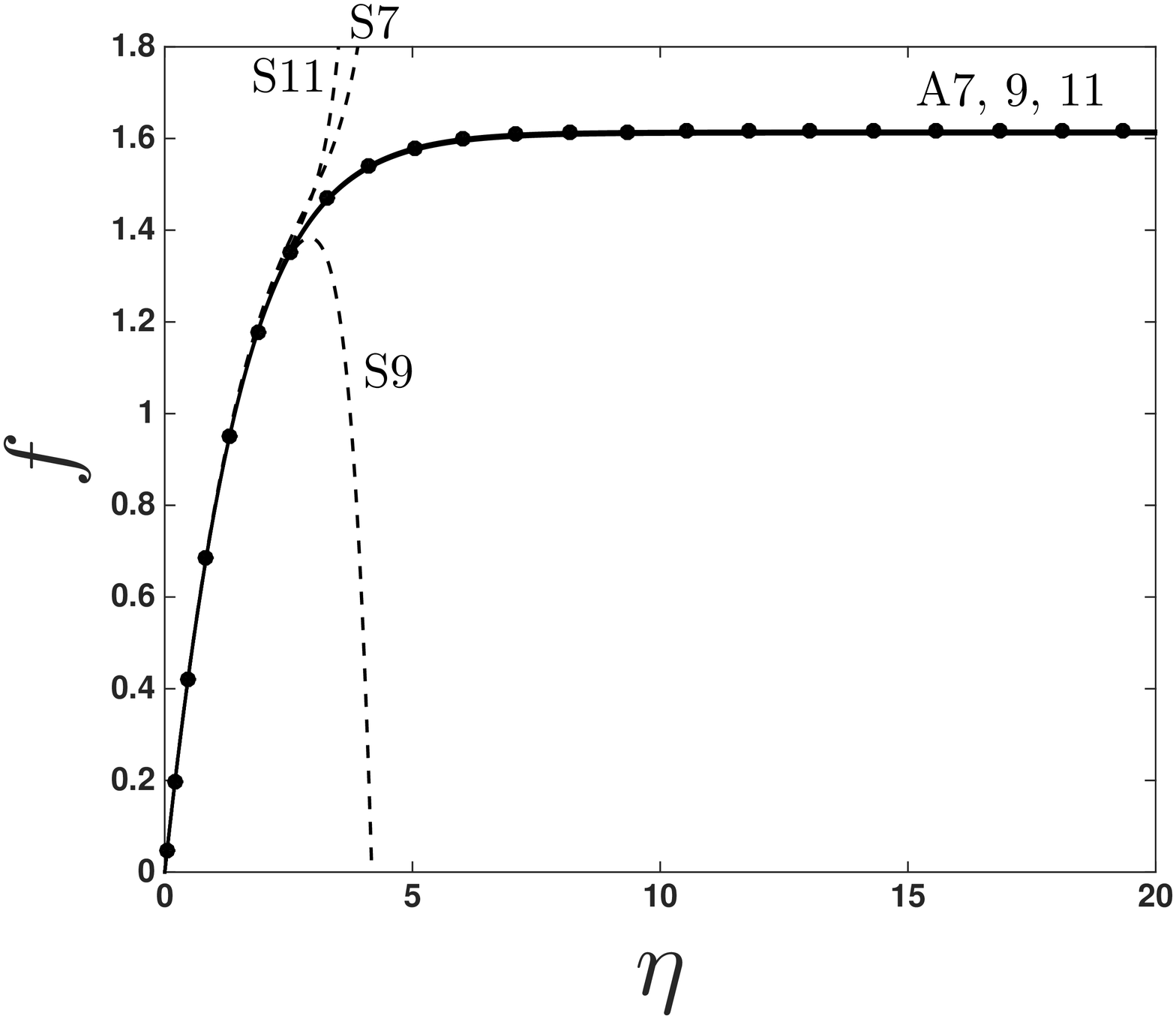}
(b)\hspace{-.05in}
\includegraphics[width=2.6in]{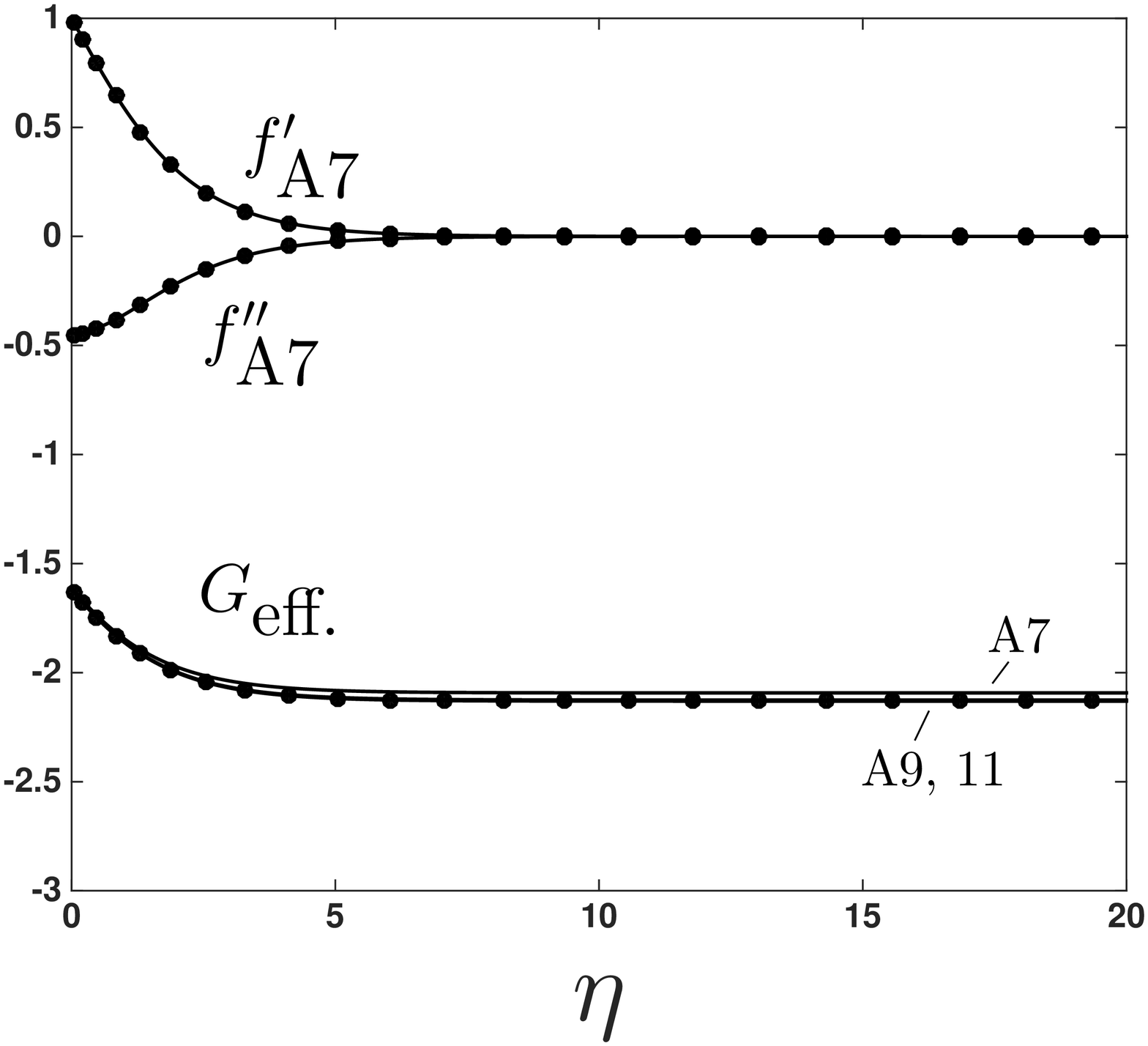}
\end{tabular}
\caption{(a) Comparison between the $N$-term series solution~(\ref{eq:series}) ($\---$) labeled as S$N$, the corresponding approximant~(\ref{Exp}) ($\--$) labeled as A$N$, and the numerical solution~\cite{Cortell} ($\bullet$) to the Sakiadis problem~(\ref{eq:Sakiadis}). (b) Derivatives of approximant~(\ref{Exp}) with $N$=7 and $G_\text{eff}$ given by~(\ref{Geff}) compared with the numerical solution~\cite{Cortell} ($\bullet$).}
\label{fig:sakiadis}
\end{figure*}

For the Sakiadis problem, a relatively simple approximant that meets the infinite condition in~(\ref{sBC}), such as approximant~(\ref{SakiadisA}), is adequate for describing $f$ and its derivatives. However, in comparing the $G_\text{eff}$ curves of Figs.~\ref{fig:sakiadisCrude}b and \ref{fig:sakiadis}b, it is clear that an approximant such as~(\ref{Exp}) which includes asymptotic corrections, is needed to predict $G$.  As we shall see in the next section, simple approximants can sometimes also be used to predict higher-order asymptotic properties, provided that the approximant converges at large enough $\eta$ for the asymptotic behavior to be reached. %- that is, before the approximant predictions of asymptotic constants diverge at larger $\eta$ (they must do so, as approximants without the correct asymptotic behavior will ultimately exhibit nonuniform convergence in the asymptotic properties.

Lastly, it is interesting to note that approximant~(\ref{Exp}) contains no singularities, and thus cannot be used to predict the closest singularity $\eta_s$ (and $S$) directly, as was done in Section~\ref{sec:Sakiadis1} for approximant~(\ref{SakiadisA}).   However, if we use either the new high precision values of $\kappa$ and $C$ from approximant~(\ref{Exp}) or the numerical values listed in the table as inputs to approximatnt~(\ref{SakiadisA}), the location of the singularity converges to the within the indicated digits to $\eta_s=-1.211393+ 3.88781$i and thus a radius of convergence of $S\equiv|\eta_s|=4.07217$.

Although the implementation is different, the approximant~(\ref{Exp}) is similar in style to that of the approximant of the  Kidder boundary value problem given recently in~\cite{IaconoBoyd2014}, where an exponential decay of the derivative at infinity is incorporated by changing the independent variable of a Pad\'e to reflect this behavior.  If convergence of this approximant for increasing order were considered, this approach would fall within the framework of the asymptotic approximants presented here. 

The errors in the approximants of this section are given in Table~\ref{table:E} to enable comparison with other approximate treatments of the Sakiadis problem; see for example~\cite{Xu}.
\begin{table}[h!]
\centering
\caption{The infinity norm of the error in the Sakiadis approximants and their 2$^\mathrm{nd}$ derivatives (to within 1 significant digit), defined respectively as $E$=$||f_A-f_\mathrm{numerical}||_\infty$ and $E_2$=$||f_A''-f_\mathrm{numerical}''||_\infty$ on the interval $0\le \eta\le20$.}
%\small
\label{table:E}
\begin{tabular}{ccccc}
\hline
$N$       & \begin{tabular}{c}$E$\\from~(\ref{SakiadisA})\end{tabular} & \begin{tabular}{c}$E_2$\\ from~(\ref{SakiadisA})\end{tabular}   & \begin{tabular}{c}$E$\\ from~(\ref{Exp})\end{tabular}   & \begin{tabular}{c}$E_2$\\ from~(\ref{Exp})\end{tabular}        \\ \hline
5         & 5$\times10^{-1}$   & 6$\times10^{-2}$ &  4$\times10^{-2}$& 2$\times10^{-2}$\\
7         & 1$\times10^{-1}$ &1$\times10^{-2}$ & 6$\times10^{-3}$& 2$\times10^{-3}$  \\
9         & 3$\times10^{-2}$ & 1$\times10^{-3}$  & 8$\times10^{-4}$&2$\times10^{-4}$\\
11        & 1$\times10^{-2}$  & 7$\times10^{-4}$  & 9$\times10^{-5}$&2$\times10^{-5}$\\
13        &  $\--$   &  $\--$   & 9$\times10^{-6}$&2$\times10^{-6}$\\
15       & $\--$    & $\--$  &8$\times10^{-7}$&2$\times10^{-7}$\\
20        & $\--$  & $\--$  &2$\times10^{-9}$&3$\times10^{-10}$\\
25        & $\--$  & $\--$   &1$\times10^{-11}$&5$\times10^{-12}$\\
30        & $\--$   & $\--$   &1$\times10^{-11}$&4$\times10^{-13}$\\
\hline
\end{tabular}
\end{table}

\subsection{The Blasius problem \label{sec:Blasius}}
The Blasius problem is the archetypal boundary-layer problem found in most undergraduate fluid mechanics books, describing the boundary layer due to a moving fluid over a stationary flat plate~\cite{Blasius}.  The differential equation in $f(\eta)$ is the same as the Sakiadis problem (and arises after a similarity transform is applied to the governing equations), 
\begin{subequations}
\begin{equation}
2f'''+ff''=0,
\label{bDE}
\end{equation}
while the boundary conditions are now
\begin{equation}
f(0)=0, ~ f'(0)=0, ~ f'(\infty) = 1.
\label{bBC}
\end{equation}
\label{eq:Blasius}
\end{subequations}
The solution to~(\ref{eq:Blasius}) may be constructed as a power series~(\ref{eq:series}) as in the Sakiadis problem. However, one difference is that the Blasius series solution skips every two terms, starting with the first two, which are zero as imposed by the boundary conditions at $\eta$=0 in~(\ref{bBC}). Again, the wall shear parameter is defined as
\[f''(0)\equiv \kappa\]
which, along with the asymptotic properties~\cite{Boyd1999}
\begin{subequations}
\label{BA}
\begin{equation}
\lim_{\eta\to\infty}(f-\eta)\equiv B
\label{BB}
\end{equation}
\begin{equation}
\lim_{\eta\to\infty}\left[\text{exp}\left(\frac{\eta^2}{4}+\frac{B\eta}{2}\right)f''\right]\equiv Q,
\label{BQ}
\end{equation}
\end{subequations}
will here be predicted by the asymptotic approximant, in a similar fashion as done by~\cite{Boyd1997} (for $\kappa$ using Pad\'es).   Recent numerical and analytical predictions of these quantities, including a review of previous results may be found in~\cite{Yun,Xu,Lal,IaconoBoyd2015}.  

The asymptotic behavior of the Blasius problem as $\eta\to\infty$ is given as (see Appendix~\ref{sec:BlasiusA}) 
\begin{equation}
f \sim \eta+B+4Q  \frac{\text{exp}[-\eta^2/4-B\eta/2]}{(\eta+B)^2}[1+O(\frac{1}{(\eta+B)^2})],~\textrm{as} \ \eta \rightarrow \infty.
\label{Q2}
\end{equation}
For comparison, the exponential correction analogous to~(\ref{Q2}) for the Sakiadis problem is~(\ref{7}) in Appendix~\ref{sec:SakiadisA}.  While we were able to incorporate~(\ref{G}) directly in the Sakiadis approximant~(\ref{Exp}), this is not an option for the Blasius problem since the asymptotic form~(\ref{Q2}) contains a singularity at $\eta=-B$ and it is known that $B<0$~\cite{Boyd1999}.  To clarify this issue, note that~(\ref{Q2}) is fully valid as $\eta\to\infty$.  However, our methodology for constructing an approximant relies on the unification of asymptotic limits such that the final form is capable of describing all $\eta>0$, and here is where~(\ref{Q2}) poses an issue in its ability to inform such an approach. 

%One may verify that~(\ref{BQ}) is consistent with~(\ref{Q2}).  

While we cannot directly use the asymptotic expansion for the Blasius problem, we can use a simple approximant that agrees asymptotically to 1$^\text{st}$ order (i.e. with (\ref{bBC}) and~(\ref{BB})) and then verify that it approaches the exponential correction~(\ref{BQ}) as $\eta\to\infty$.  We use an approximant of the form
\begin{subequations}
 \label{BlasiusA}
 \begin{equation}
  f_A=\eta+B-B\left(1+\sum_{n=1}^{N}A_n~\eta^n\right)^{-1},
 \label{BlasiusApproximant}
 \end{equation}
which automatically satisfies both~(\ref{bBC}) and~(\ref{BB}).  As similarly discussed in the context of the Sakiadis Approximant~(\ref{SakiadisApproximant}), the above form is not a Pad\'e approximant.  Firstly, $B$ is an unknown parameter in the approximant, and so~(\ref{BlasiusApproximant}) cannot be considered a Pad\'e for $f-\eta-B$. Secondly, if one combines the terms of~(\ref{BlasiusApproximant}) through a common denominator, the coefficients of the numerator will have an explicit dependence on those in the denominator, whereas this is not the case for standard Pad\'es.  

The coefficients $A_0\dots A_{N}$, $B$, and $\kappa$ are now calculated such that the $N$-term Taylor expansion of~(\ref{BlasiusApproximant}) about $\eta=0$ is exactly equal to the $N$-term truncation of~(\ref{eq:series}). Following the same procedure as in Section~\ref{sec:Sakiadis1}, we arrive at the following recursion for the coefficients, which are now functions of $\kappa$ and $B$:
\begin{equation}
A_{n>0}=\frac{1}{B}\sum_{j=1}^n\tilde{a}_{j}~A_{n-j},~~A_0=1,
\label{BlasiusCoefficients}
\end{equation}
where $\tilde{a}_1=-1$ and $\tilde{a}_{j>1}=a_j$.  The coefficients $A_N$ and $A_{N-1}$ are now sacrificed to simultaneously predict $\kappa$ and $B$.  Setting  
\begin{equation}
A_N(\kappa,B)=0,~A_{N-1}(\kappa,B)=0
\end{equation}
\end{subequations}
in~(\ref{BlasiusCoefficients}) leads to two nonlinear equations, whose most rapidly converging $\kappa$ and $B$ roots are given in Table~\ref{table:Bkappa}, along with the magnitude of the singularity closest to $\eta=0$ in the approximant, denoted by $S$.  The table reports values up to $N=50$, beyond which convergence cannot be established in the 8th digit.  Our predictions converge to within 7 digits of the benchmark value for $\kappa$, 6 digits for $B$, and 2 digits for the radius of convergence $S$.  Again, we determine $S$ as the magnitude of the closest singularity, $\eta_s$, from the origin in the complex $\eta$ plane.  Although not shown in the table, the approximant~(\ref{BlasiusA}) is consistent with the literature~\cite{Boyd1999} in that, as $N$ is increased, it predicts $\eta_s$ to lie on the negative real axis.

\begin{table}[h!]
\centering
\caption{Predictions from the Blasius approximant~(\ref{BlasiusA}). The infinity norm of the error in the approximant and its 2$^\mathrm{nd}$ derivative are recorded below (to within 1 significant digit), defined respectively as $E$=$||f_A-f_\mathrm{numerical}||_\infty$ and $E_2$=$||f_A''-f_\mathrm{numerical}''||_\infty$ on the interval $0\le \eta\le8.8$.}  %The first numerical $R$ value is computed by the authors, using the numerical $\kappa$ and applying the ratio test to~(\ref{eq:series}) and the 2nd $R$ value is from~\cite{Boyd2008}.}
%\small
\label{table:Bkappa}
\setlength{\tabcolsep}{.4em}
\small
\begin{tabular}{cccccc}
\hline
$N$         & $\kappa$          & $B$              & $S$         & $E$ & $E_2$\\  \hline 
5         & 0.12545065        & -6.08950239      & 9.85302    &3$\times10^{0}$ & 2$\times10^{-1}$ \\
10        & 0.30018461        & -2.00003122      & 4.47053    &3$\times10^{-1}$& 3$\times10^{-2}$  \\
15        & 0.32626056        & -1.77168204      & 4.59756    &5$\times10^{-2}$& 5$\times10^{-3}$  \\
20        & 0.33090243        & -1.73173129      & 4.92273     &1$\times10^{-2}$ & 1$\times10^{-3}$\\
25        & 0.33181978        & -1.72323142      & 5.18569     &2$\times10^{-3}$&2$\times10^{-4}$ \\
30        & 0.33200793        & -1.72133641      & 5.43362     &5$\times10^{-4}$ &5$\times10^{-5}$ \\
35        & 0.33204717        & -1.72090879      & 5.64460     &1$\times10^{-4}$ &1$\times10^{-5}$ \\
40        & 0.33205518        & -1.72081494      & 5.69196       &3$\times10^{-5}$ &2$\times10^{-6}$\\
45        & 0.33205693        & -1.72079309      & 5.68929      &5$\times10^{-6}$ &4$\times10^{-7}$ \\
50        & 0.33205731        & -1.72078801      & 5.68933     &3$\times10^{-7}$ &3$\times10^{-8}$ \\
numerical~\cite{Lal,Boyd2008} & 0.332057336215196 & -1.7207876575205 & 5.6900380545 &0&0\\ \hline
\end{tabular}
\end{table}

Once $B$ and $\kappa$ are substituted into~(\ref{BlasiusCoefficients}) for all non-zero coefficients ($A_1$ through $A_{N-2}$), the approximant series in~(\ref{BlasiusA}) may be constructed.  The approximant~(\ref{BlasiusA}) is compared with the series and numerical solution of~(\ref{eq:Blasius}) in Fig.~\ref{fig:blasius}.    Like the simple Sakiadis approximant given in Section~(\ref{sec:Sakiadis1}), this basic Blasius approximant converges surprisingly well, as can be seen in Fig.~\ref{fig:blasius}a.  In Fig.~\ref{fig:blasius}b, we confirm that $f'$ and $f''$ are accurately obtained from the approximant, once $f$ has reasonably converged (here, for $N=30$).  Unlike the simple Sakiadis approximant, the Blasius approximant is capable of picking up a higher-order asymptotic quantity - namely $Q$ defined in~(\ref{BQ}).  This is seen in Fig.~\ref{fig:blasius}c, where an effective $Q$, defined as 
 \begin{equation}
 Q_\text{eff}=\text{exp}\left(\frac{\eta^2}{4}+\frac{B\eta}{2}\right)f'',~Q_\text{eff}\to Q~\text{as}~\eta\to\infty
 \label{Qeff}
 \end{equation}
is plotted versus $\eta$ for various series and approximant order.  While the approximant converges directly to a constant $Q$ value at large $\eta$, it is striking that the series solution also reaches the $Q$ plateau prior to diverging.  This may be a result of the Blasius series having a region of convergence that overlaps with asymptotic effects felt at relatively small $\eta$.  As seen in Fig.~\ref{fig:blasius}c, both the approximant and the numerical solution converge to $Q\approx0.1115$, which is consistent with the value $Q$=0.111483755 obtained numerically by~\cite{Lal} and semi-analytically (i.e.~dependent on numerical $\kappa$) in~\cite{IaconoBoyd2015}\footnote{$Q$ in~\cite{IaconoBoyd2015} is defined as $e^{-B^2/4}$ multiplied by the $Q$ defined here.}; in both references, the value is reported well-beyond these digits, but both round to the value above.  Since the behavior~(\ref{Qeff}) is not explicitly incorporated into the approximant~(\ref{BlasiusA}), the convergence of $Q_\text{eff}$ is not uniform at large $\eta$.  For example, in Fig.~\ref{fig:blasius}c, the behavior of the $N=20, 30, 40$ approximants (labeled A20, A30, A40) demonstrates the non-uniform convergence that will occur for any order $N$ at some value of $\eta$.  

\begin{figure*}[h!]
\begin{tabular}{cc}(a)\hspace{-.05in}
\includegraphics[width=2.6in]{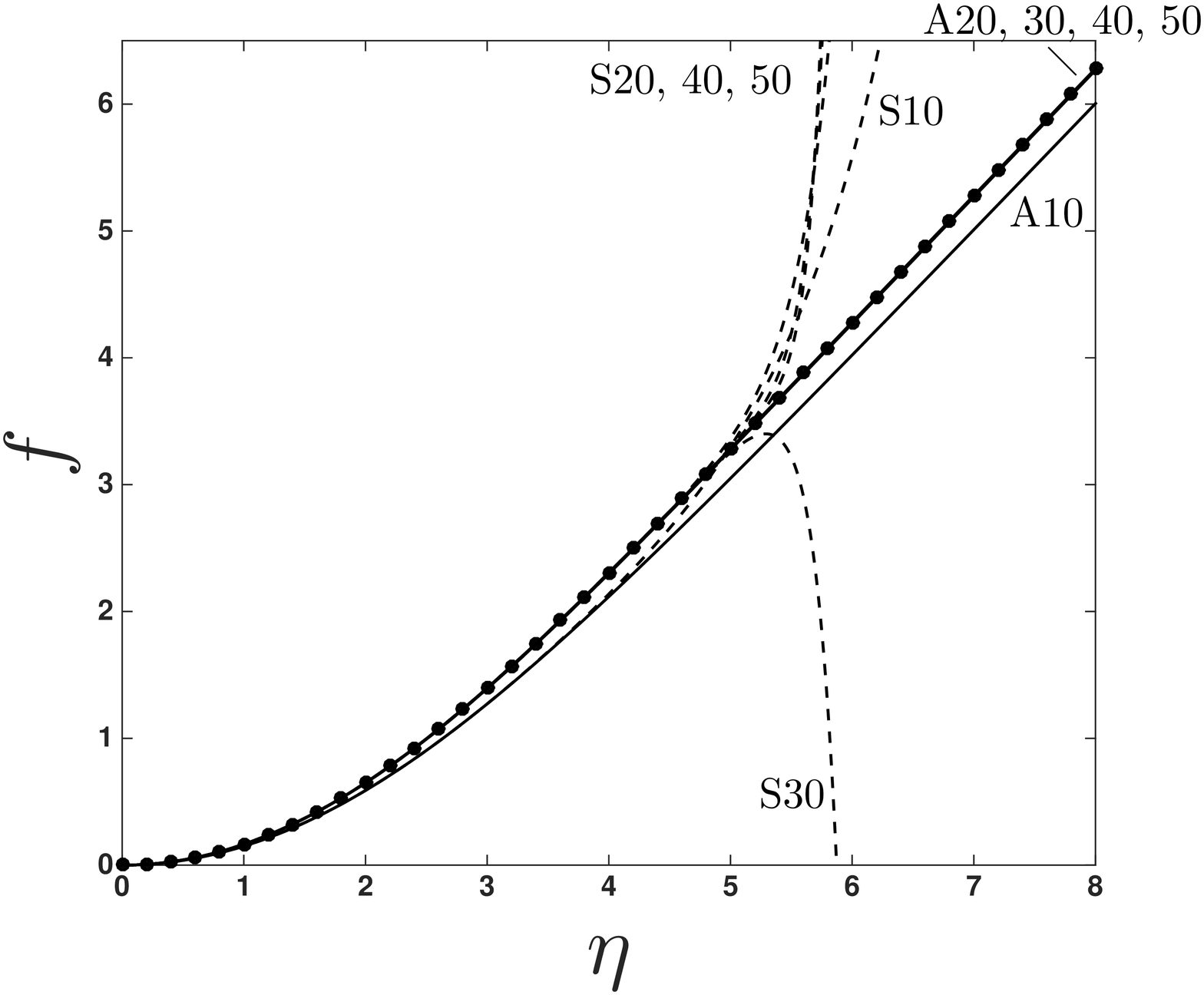}
(b)\hspace{-.05in}
\includegraphics[width=2.6in]{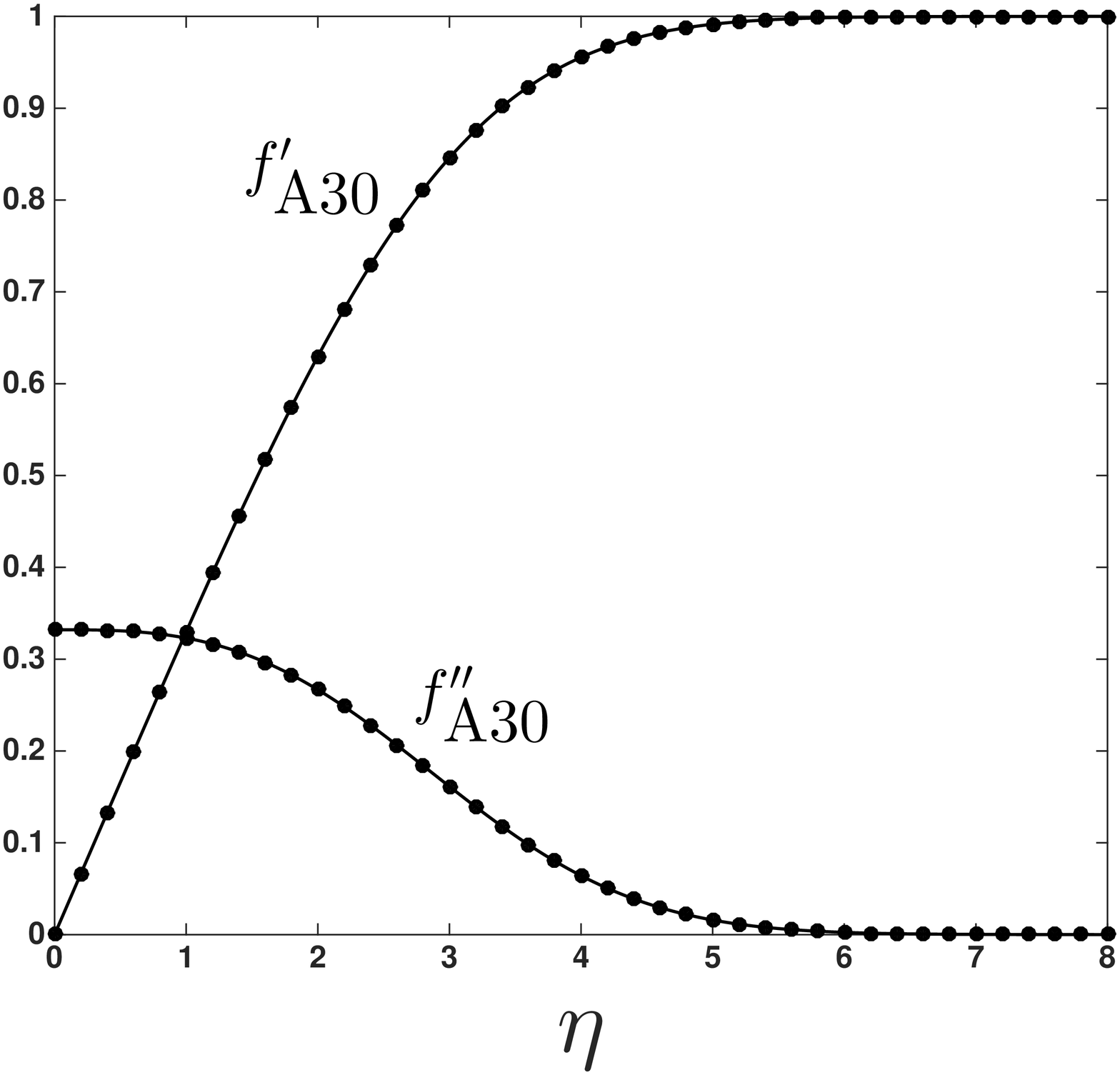}\\
(c)\includegraphics[width=2.6in]{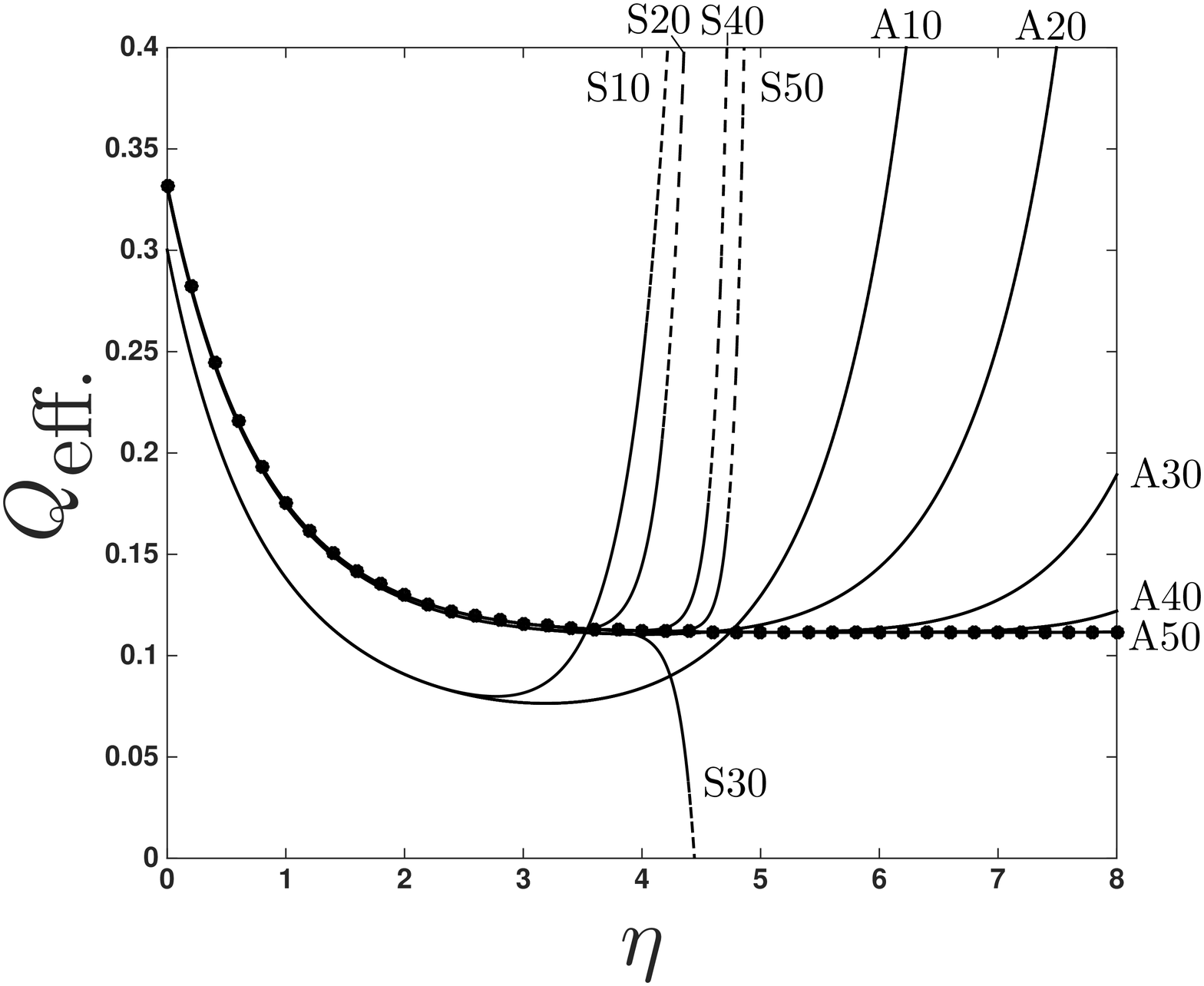}
\end{tabular}
\caption{(a) Comparison between the $N$-term series solution~(\ref{eq:series}) ($\---$) labeled as S$N$ and the corresponding approximant~(\ref{BlasiusA}) ($\--$) labeled as A$N$ to the Blasius problem~(\ref{eq:Blasius}). (b) Derivatives of approximant~(\ref{BlasiusA}) with $N$=30. (c) $Q_\text{eff}$ given by~(\ref{Qeff}) as predicted by the series~(\ref{eq:series})  and approximant~(\ref{BlasiusA}). For comparison, the numerical solution from~\cite{Lal} is shown in all above figures as $\bullet$s. }
\label{fig:blasius}
\end{figure*}

The errors in the approximant are given in the last two columns of Table~\ref{table:Bkappa} to enable comparison with other approximate treatments of the Blasius problem~\cite{Boyd1997,Boyd1999,Yun,Xu,IaconoBoyd2015}. 

The virtue of approximant~(\ref{BlasiusA}) is that it is globally accurate from small to large $\eta$, is able to predict $\kappa$ and $B$ (i.e. does not rely on these values as inputs), and has coefficients that are generated from a simple recursion.  A recent alternative globally accurate approximation to the Blasius solution is provided in~\cite{IaconoBoyd2015}.  An attractive feature of this alternative approximation is that it has a simpler form than the asymptotic approximant~(\ref{BlasiusA}).  However, the form in~\cite{IaconoBoyd2015} requires that $\kappa$ is known beforehand.  Nevertheless, one may directly compare the error in $f_A''$ given in Table~\ref{table:Bkappa} with that given in figure 2 of~\cite{IaconoBoyd2015}.  

%As one might expect, as the simplicity of our approximant decreases (for increasing $N$), the accuracy increases and eventually outperforms the approximations of~\cite{IaconoBoyd2015}.  

\subsection{The Flierl-Petviashvili problem \label{sec:FP}}
The Flierl-Petviashvili (FP) equation in similarity variables $u(r)$ is used to describe vortex solitons in the ocean, atmosphere~\cite{Flierl}, and Jupiter's red spot~\cite{Petviashvili}, and is given by 
\begin{subequations}
\begin{equation}
u''+\frac{1}{r}u'-u-u^2=0
\end{equation}
\begin{equation}
u'(0)=0,~~u(\infty)=0.
\label{FPBC}
\end{equation}
\label{FP}
\end{subequations}
The coefficients of the power series solution 
\begin{subequations}
\label{FPSeries}
\begin{equation}
u= \sum\limits_{n=0}^\infty a_n r^n
\end{equation}
 to~(\ref{FP}) are given by
\begin{equation}
a_{n+2}=\frac{a_n+\sum\limits_{k=0}^n a_ka_{n-k}}{(n+2)^2},
\label{FPseries}
\end{equation}
\end{subequations}
which requires the specification of the first two coefficients to generate the remaining even coefficients; all remaining odd terms are zero.  The coefficient $a_1$=0 is known from the boundary conditions~(\ref{FPBC}) and $a_0$ is an unknown parameter,
\[a_0\equiv z\]
which, along with the asymptotic property
\begin{equation}
\lim_{r\to\infty}(ue^r\sqrt{r})\equiv D
\label{D}
\end{equation}
will be predicted by the approximant, in a similar fashion as done by~\cite{Boyd1997,Wazwaz} (for $z$ using Pad\'es).  The property~(\ref{D}) above follows from an asymptotic expansion of the FP problem as $r\to\infty$ (derived in Appendix~\ref{sec:FPA}):
\begin{equation}
u \sim D \frac{e^{-r}}{\sqrt{r}} [1+O(\frac{1}{r})],~\textrm{ as } r \rightarrow \infty.
\label{FPasymp}
\end{equation}
As was the case for the Blasius problem, a singularity in the asymptotic form (here at $r$=0 in the above) prevents direct incorporation of~(\ref{FPasymp}) to construct an approximant.  However, like the Blasius approximant of Section~\ref{sec:Blasius}, we shall again use a simple approximant that agrees asymptotically to zeroth order (i.e. with the 2$^\mathrm{nd}$ condition of (\ref{FPBC})) and then verify that it approaches the  correction~(\ref{FPasymp}) as $r\to\infty$. 

We consider a simple approximant of the form 
 \begin{subequations}
 \label{FPA}
 \begin{equation}
  u_A=\frac{z}{1+\sum\limits_{n=1}^{N}A_n~r^n},
 \label{FPApproximant}
 \end{equation}
which automatically satisfies both boundary conditions of~(\ref{FP}), including the condition $u(\infty)=0$ not explicitly captured by~(\ref{FPSeries}).  Note that~(\ref{FPApproximant}) is a Pad\'e.  The simple structure of this Pad\'e, however, allows for the unknowns to be calculated without inverting a matrix.  The coefficients $A_0\dots A_{N}$ and $z$ in~(\ref{FPApproximant}) are now calculated such that the $N$-term Taylor expansion of~(\ref{FPApproximant}) about $r=0$ is exactly equal to the $N$-term truncation of~(\ref{FPSeries}). Following the same procedure as in Section~\ref{sec:Sakiadis1}, we arrive at the following recursion for the coefficients, which are now functions of $z$:
\begin{equation}
A_{n>0}=-\frac{1}{z}\sum_{j=1}^na_{j}~A_{n-j},~~A_0=1.
\label{FPCoefficients}
\end{equation}
The coefficients above mimic the series~(\ref{FPSeries}) in that it skips odd coefficients. The unknown $z$ in~(\ref{FPCoefficients}) is determined by sacrificing the degree of freedom normally used to compute $A_N$. This is equivalent to setting $A_N=0$ in~(\ref{FPCoefficients}), leading to 
\begin{equation}
\sum_{j=1}^Na_jA_{N-j}=0.
\label{zeqn}
\end{equation}
\end{subequations}
Note that since the $a_j$ coefficients are polynomials in $z$, (\ref{zeqn}) is also a polynomial in $z$.  The most rapidly converging $z$-roots of~(\ref{zeqn}) (for increasing $N$) are given in Table~\ref{table:z}.  Once $z$ is substituted into~(\ref{FPCoefficients}) for all non-zero coefficients ($A_1$ through $A_{N-1}$), the approximant series in~(\ref{FPApproximant}) may be constructed.  The magnitude of the singularity  closest to $r=0$ in the approximant~(\ref{FPA}) is also given in Table~\ref{table:z}, denoted by $S$.   Although not shown in the table, the approximant~(\ref{FPApproximant}) predicts that, within the reported digits, the singularity appears to lie on the positive imaginary axis.

\begin{table}[h!]
\centering
\caption{Predictions from the Flierl-Petviashvili approximant~(\ref{FPApproximant}).  The numerical method for computing the values in the bottom row is described in Appendix~\ref{sec:FPnum}. The infinity norm of the error in the approximant is also recorded below (to within 1 significant digit), defined as $E$=$||u_A-u_\mathrm{numerical}||_\infty$ on the interval $0\le r\le 10$. }
%\small
\label{table:z}
\begin{tabular}{cccc}
\hline
$N$         & $z$    & $S$ & $E$\\ \hline
%2         & -1    \\
4         & -1.5   & 2.82843 &  9$\times10^{-1}$\\
6         & -2.14039 & 2.57284 & 3$\times10^{-1}$  \\
8         & -2.34792  & 2.47227  & 4$\times10^{-2}$\\
10        & -2.38564  & 2.60878  & 6$\times10^{-3}$\\
12        & -2.39117  & 2.59927  & 7$\times10^{-4}$\\
14        & -2.39185   & 2.61160 &2$\times10^{-4}$\\
16        & -2.39196  & 2.61154  &1$\times10^{-5}$\\
%18        & -2.39193    \\
numerical & -2.3919564032 & 2.611541077& 0 \\ \hline
\end{tabular}
\end{table}

The approximant~(\ref{FPA}) is compared with the series and numerical solution of~(\ref{FP}) in Fig.~\ref{fig:FP}.   As for the simple Sakiadis approximant given in Section~\ref{sec:Sakiadis1} and the Blasius approximant of Section~\ref{sec:Blasius}, the FP approximant also converges uniformly, as can be seen in Fig.~\ref{fig:FP}a.  In Fig.~\ref{fig:FP}b, we confirm that $u'$ is accurately obtained from the approximant once $u$ has reasonably converged (here, for $N=12$).  Like the Blasius approximant, the FP approximant is capable of picking up a higher-order asymptotic quantity - here,  $D$ defined in~(\ref{D}).  This is seen in Fig.~\ref{fig:FP}b, where an effective $D$, defined as 
 \begin{equation}
 D_\text{eff}=ue^r\sqrt{r},~D_\text{eff}\to D~\text{as}~r\to\infty
 \label{Deff}
 \end{equation}
 is plotted versus $r$ for various approximant order.  As $N$ is increased, the approximant $D_\text{eff}$ appears to be converging to the numerical solution, which is in turn converging to a value of
 \[D\approx-10.7\]  
which may serve as a useful metric for future analytical and numerical solutions of the FP problem. Again, Fig.~\ref{fig:FP} shows that convergence in $D_\text{eff}$ is not uniform at large $r$, as expected by the assumed form of the approximant.

\begin{figure*}[h!]
\begin{tabular}{cc}(a)\hspace{-.05in}
\includegraphics[width=2.6in]{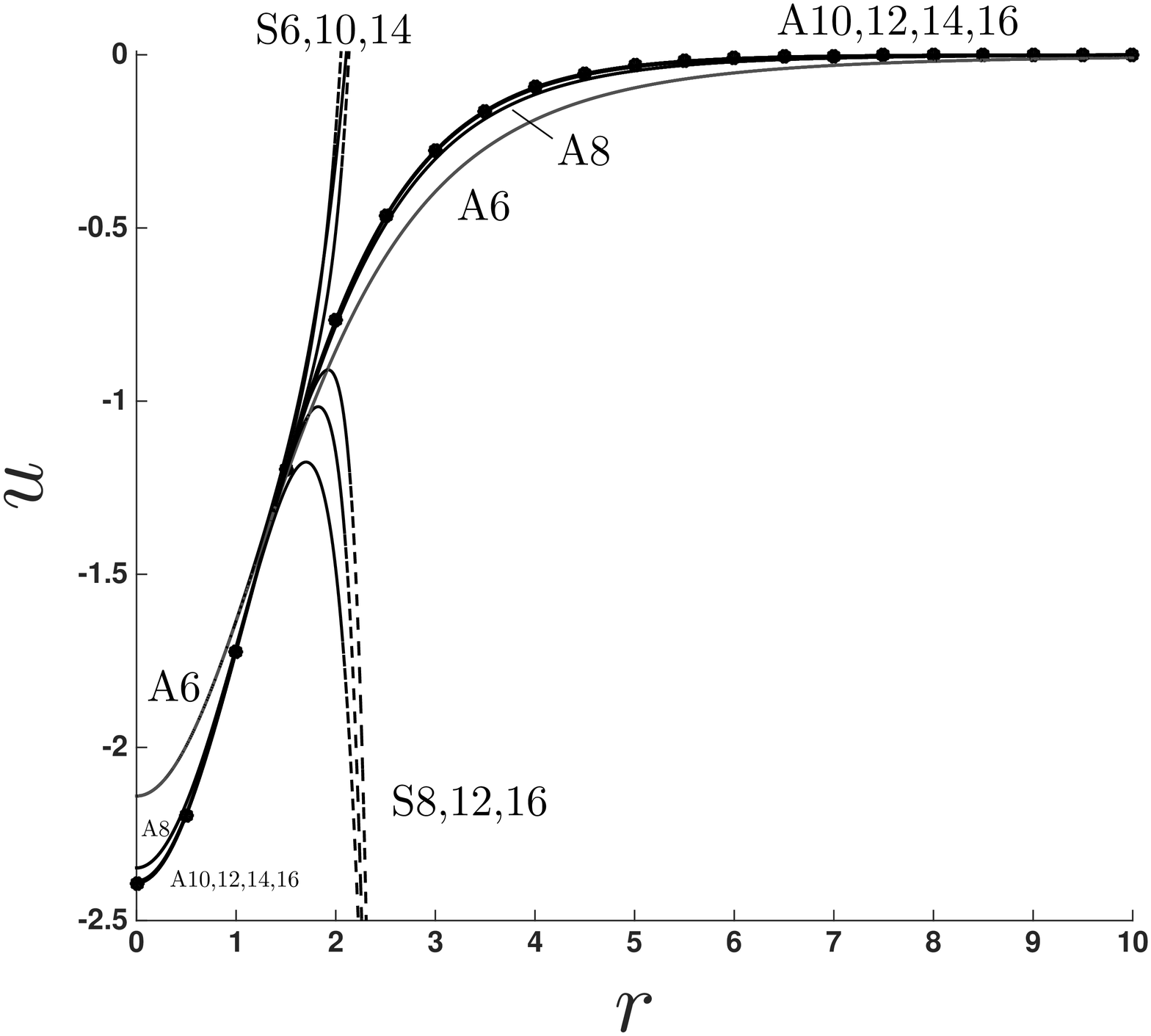}
(b)\hspace{-.05in}
\includegraphics[width=2.6in]{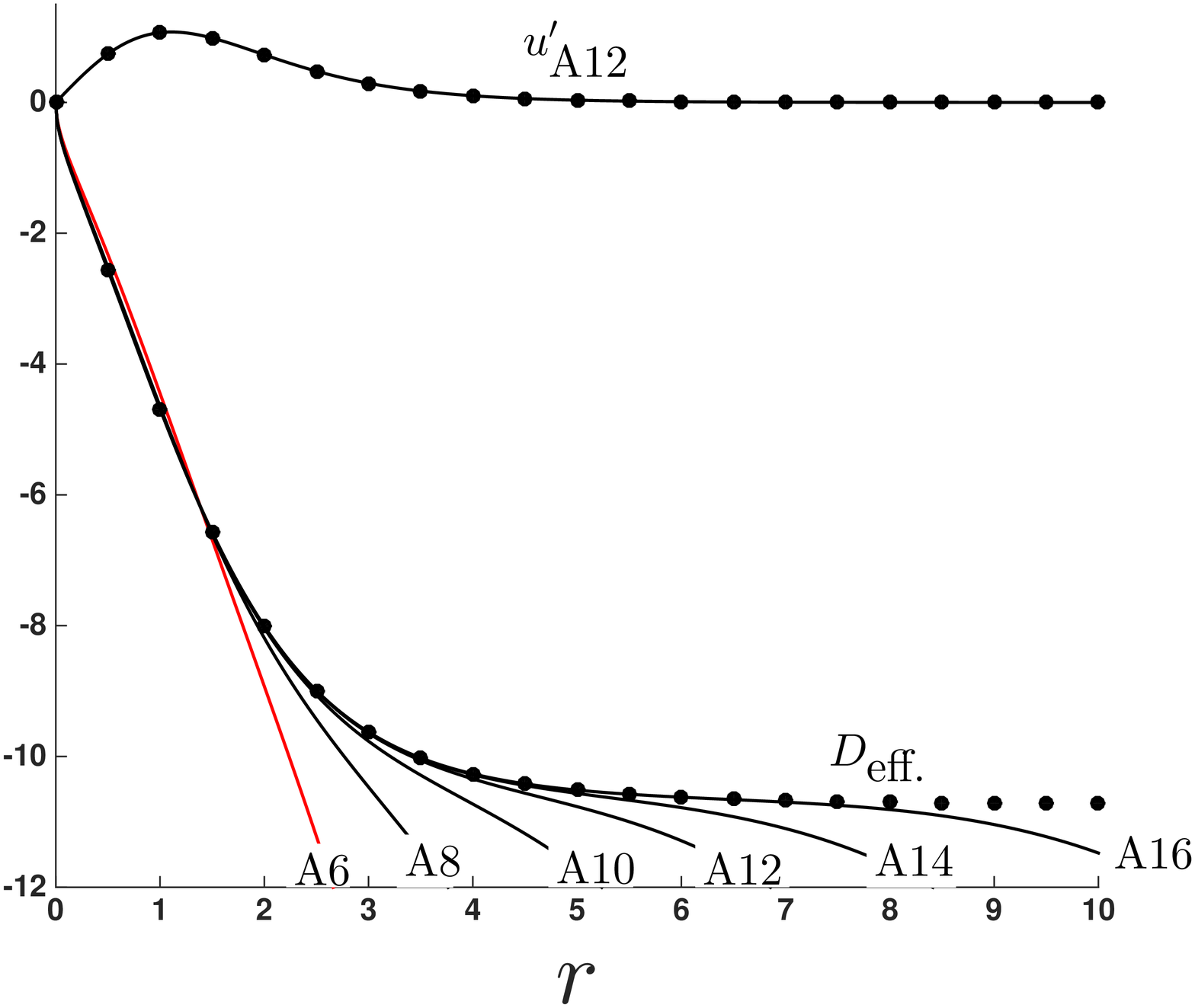}
\end{tabular}
\caption{(a) Comparison between the $N$-term series solution~(\ref{FPSeries}) ($\---$) labeled as S$N$ and the corresponding approximant~(\ref{FPA}) ($\--$) labeled as A$N$ to the Flierl-Petviashvilli problem~(\ref{FP}). (b) Derivatives of approximant~(\ref{FPA}) with $N$=12 and $D_\text{eff}$ given by~(\ref{Deff}). For comparison, the numerical solution (described in Appendix~\ref{sec:FPnum}) is shown in all above figures as $\bullet$s. }
\label{fig:FP}
\end{figure*}

The errors in the approximant are given in the rightmost column of Table~\ref{table:z} to enable comparison with other approximate treatments of the FP equation; see for example~\cite{Boyd1997}.

\section{Summary \label{sec:summary}}

In this work, we formalize a new approach to sum series.   We provide both a methodology and examples that demonstrate how asymptotic approximants may be constructed.  The key feature is that such approximants are designed so that the asymptotic behaviors in two regions of the domain may be joined. This approach has seen recent success in the analytic continuation of ``virial'' series that describe the pressure-density-temperature dependence of various model fluids, where no differential equation is available and usually only the first 3 to 12 terms of these (often divergent) power series are known~\cite{BarlowJCP, BarlowAIChE, Barlow2015}. 

 Here, we provide an additional application for asymptotically consistent approximants - namely, a remedy for divergent power series that arise as solutions to nonlinear ordinary differential equations. The analytic forms provided here enable symbolic differentiation and thus allow analytical evaluation of the flow field (shear stress, vorticity, etc.) to any desired resolution at low computational expense.  Additionally, the approximant coefficients are relatively straightforward to obtain, as they often only require the implementation of a few series identities, enabling them to be generated typically via a simple recurrence relation.  
 
Note that improvements to the forms given here are easily made, depending on how many parameters are to be embedded in (or predicted by) the approximant.  Here, we use the approximant to not only find a solution, but also provide the as-of-yet most accurate and precise values of the wall shear, asymptotic constants, and singularity of smallest magnitude for the Sakiadis boundary layer.  If better estimates become available from numerical (or other) techniques, these may be used as inputs to the approximant to improve global accuracy.  For the Blasius problem and Flierl-Petviashvilli problem, more accurate values for some of these properties are, in fact, available and could have been used as inputs to the approximants; we chose instead to predict these quantities to demonstrate that the method is self contained and does not rely on external data.  

%Nevertheless, the approximants in this paper predict accurate values for these two problems, confirming recent numerical studies.   

The results presented here motivate the development and application of asymptotic approximants to other problems where divergent, truncated, and underspecified series arise.  Considering the capabilities afforded and the ease with which asymptotic approximants are generated, they deserve consideration as an alternative to the more traditional Pad\'e approximants in problems of mathematical physics.

% On the other hand, one may wish to swap an unknown coefficient of the approximant to instead predict $\kappa$ for either the Blasius or Sakiadis problems.  Here a mixed approach was taken, since accurate values of $\kappa$ are readily available.  

\bibliographystyle{unsrt}
\bibliography{Sakiadis}

\appendix
\section{Useful Series Formulae\label{sec:formulae}}
The following relations may be used to develop a recursion for the coefficients of an asymptotic approximant, allowing one to avoid solving an algebraic system.  The first relation is the well-known Cauchy product of two series~\cite{Churchill}:
\begin{equation}
\sum_{n=0}^Na_nx^n\sum_{n=0}^NA_nx^n=\sum_{n=0}^N\left(\sum_{j=0}^na_jA_{n-j}\right)x^n. 
\label{Cauchy}
\end{equation}
Setting both sides of~(\ref{Cauchy}) equal to one and rearranging, the recursion leads to a representation for the $N$-term expansion of the reciprocal of a series:
\begin{subequations}
\label{inverse}
\begin{equation}
\left(\sum_{n=0}^Na_nx^n\right)^{-1}=\sum_{n=0}^NA_nx^n,
\label{inverseseries}
\end{equation}
where 
\begin{equation}
A_{n>0}=-\frac{1}{a_0}\sum_{j=1}^na_jA_{n-j},~A_0=\frac{1}{a_0}.
\end{equation}
\end{subequations}
The generalization of~(\ref{inverse}) for the $N$-term expansion of a series raised to any real power $s$ is given by J. C. P. Miller's formula~\cite{Henrici}:
\begin{subequations}
\label{Miller}
\begin{equation}
\left(\sum_{n=0}^Na_nx^n\right)^{s}=\sum_{n=0}^NA_nx^n,
\end{equation}
where 
\begin{equation}
A_{n>0}=\frac{1}{n~a_0}\sum_{j=1}^n(js-n+j)a_jA_{n-j},~A_0=(a_0)^s.
\end{equation}
\end{subequations}

\section{ Asymptotic Expansions \label{sec:asymptotics}}
\subsection{Sakiadis Problem \label{sec:SakiadisA}}
We wish to determine the $\eta \rightarrow \infty$ asymptotic behavior of the solution $f(\eta)$ to the Sakiadis equation~(\ref{eq:Sakiadis}).  For purposes of the analysis, the infinite constraint~(\ref{C}) alone suffices. We write:
\begin{equation}
f = C + h(\eta),\ \textrm{with} \ h \rightarrow 0\ \textrm{as} \ \eta \rightarrow \infty
\label{1-2}
\end{equation}
where $C$ is a constant and $h(\eta)$ is a function to be determined. The form~(\ref{1-2}) is substituted into equation~(\ref{sDE}) to obtain
\begin{equation}
2h'''+(C+h)h'' = 0.
\label{B2}
\end{equation}
The equation~(\ref{B2}) may be simplified by noting that terms quadratic in $h$ are small relative to other terms as $\eta\to\infty$, i.e.:
\begin{equation}
2h'''+ Ch'' = 0,~ \  \textrm{as} \ \eta \rightarrow \infty.
\label{heqn}
\end{equation}
The solution of this equation is
\begin{equation}
h \sim G e^{-C\eta/2},~\ \textrm{as} \ \eta \rightarrow \infty,
\end{equation}
where $G$ is a constant. Higher order corrections can be obtained by assuming the form
\begin{equation}
h \sim G e^{-C\eta/2} + D(\eta),~\ \textrm{as} \ \eta \rightarrow \infty,
\label{1-6}
\end{equation}
and applying the method of dominant balance~\cite{Bender} when substituted into~(\ref{B2}), which leads to 
\begin{equation}
D(\eta)\sim \frac{G^2}{4C} e^{-C\eta},~\ \textrm{as} \ \eta \rightarrow \infty.
\label{1-8}
\end{equation}
The asymptotic behavior is written concisely by combining equation~(\ref{1-2}),~(\ref{1-6}), and~(\ref{1-8}) to yield:
\begin{equation}
f \sim C+ G e^{-C\eta/2} + \frac{G^2}{4C} e^{-C\eta} + O(e^{-3C\eta/2}),\ \textrm{as} \ \eta \rightarrow \infty.
\label{7}
\end{equation}
The above process for obtaining corrections may be repeated and leads to a series of exponentials with arguments $-nC\eta/2$ for $n=1, 2, 3,4 \dots$.

\subsection{Blasius Problem \label{sec:BlasiusA}}
We wish to determine the $\eta \rightarrow \infty$ asymptotic behavior of the solution $f(\eta)$ to the Blasius problem~(\ref{eq:Blasius}).  For purposes of the analysis, the infinite constraint~(\ref{BB}) alone suffices. We write:
\begin{equation}
f = \eta + B + g(\eta),\ \textrm{with} \ g \rightarrow 0\ \textrm{as}\ \eta \rightarrow \infty,
\label{2-2}
\end{equation}
where \textit{B} is an unknown constant and $g(\eta)$ is a function to be determined. The form~(\ref{2-2}) is substituted into~(\ref{eq:Blasius}) to obtain
\begin{equation}
2g'''+(\eta+B+g)g'' = 0.
\label{B9}
\end{equation}
The equation~(\ref{B9}) may be simplified by noting that terms quadratic in $g$ are small relative to other terms as $\eta\to\infty$, i.e.:
\begin{equation}
2g'''+(\eta+B)g'' = 0, \ \textrm{as}\ \eta \rightarrow \infty.
\end{equation}
The solution of this equation is
\begin{equation}
g''\sim Q e^{[-(\eta^2 + 2B\eta)/4]},~\textrm{as} \ \eta \rightarrow \infty,
\end{equation}
where \textit{Q} is an unknown constant.  We can then write:
\begin{equation}
g'\sim Q e^{(B^2/4)}\int_{\infty}^{\eta} e^{[-(\xi+B^2)/4]}d\xi ,~\textrm{as} \ \eta \rightarrow \infty
\end{equation}
and after application of integration by parts, we obtain:
\begin{equation}
g' \sim -2Q e^{(B^2/4)} \frac{e^{[-(\eta+B^2)/4]}}{\eta+B}[1+O(\frac{1}{(\eta+B)^2})],~\textrm{as} \ \eta \rightarrow \infty.
\end{equation}
The process can then be repeated to find an expression for $g$, i.e.:
\begin{equation}
g \sim2Q e^{(B^2/4)} \int_{\infty}^{\eta}\frac{e^{[-(\xi+B^2)/4]}}{\xi+B}d\xi,~\textrm{as} \ \eta \rightarrow \infty,
\end{equation}
from which integration by parts yields
\begin{equation}
g \sim4Q e^{(B^2/4)} \frac{e^{[-(\eta+B^2)/4]}}{(\eta+B)^2}[1+O(\frac{1}{(\eta+B)^2})],~ \textrm{as} \ \eta \rightarrow \infty,
\end{equation}
which can be combined with the original assumed form~(\ref{2-2}) to obtain the asymptotic behavior:
\begin{equation}
f \sim \eta+B+4Q  \frac{\exp[-\eta^2/4-B\eta/2]}{(\eta+B)^2}[1+O(\frac{1}{(\eta+B)^2})],~\textrm{as} \ \eta \rightarrow \infty.
\end{equation}

\subsection{Flierl-Petviashvili Problem \label{sec:FPA}}
We wish to determine the $r\rightarrow \infty$ asymptotic behavior of the solution $u(r)$ to the Flierl-Petviashvili problem~(\ref{FP}).  For purposes of the analysis, the infinite constraint in~(\ref{FPBC}) alone suffices. Since $u\rightarrow 0$ as $r\rightarrow \infty$, it is seen by inspection of~(\ref{FP}) that $u^2 <<u$ as $r\rightarrow \infty$, leading to the linear equation
\begin{equation}
r^2 u'' + ru - r^2u = 0, ~\textrm{ as } r \rightarrow \infty.
\label{3-3}
\end{equation}
The solution of~(\ref{3-3}) is a modified Bessel function of zeroth order~\cite{Abramowitz}(p. 374 equation 9.6.1). The leading-order behavior of this function is given by equation 9.7.2 (p. 378) of~\cite{Abramowitz} as 
\begin{equation}
u \sim D \frac{e^{-r}}{\sqrt{r}} [1+O(\frac{1}{r})],~\textrm{ as } r \rightarrow \infty,
\end{equation}
where $D$ is an unknown constant. 
%\numberwithin{equation}{section}
% \setcounter{equation}{0}

\section{Numerical Solution of the Flierl-Petviashvili Problem \label{sec:FPnum}}
The numerical solution to the Flierl-Petviashvili equation used in Table~\ref{table:z} is obtained by first applying the transformation $\alpha = 1-e^{-r}$ to~(\ref{FP}), leading to 
\begin{eqnarray}
\nonumber
(1-\alpha)^2 \frac{d^2u}{d\alpha ^2}-[(1-\alpha)+\frac{1+\alpha}{\ln (1-\alpha)}]\frac{du}{d\alpha}-u-u^2=0\\
\frac{du}{d\alpha}(0)=0,~u(1)=0.
\label{FP2}
\end{eqnarray}
Note that, like the original FP equation~(\ref{FP}), the solution to~(\ref{FP2}) is not unique. There is a trivial solution $u=0$ and a non-trivial solution - we seek the latter.  The ``Shooting Method''~\cite{Turner} (i.e. iterating on $u(0)=z$ guesses) is used to effectively convert~(\ref{FP2}) into an initial value problem, which is solved via 4th order Runge-Kutta.  The boundary condition $u(1)=0$, which informs the shooting iteration, is replaced with $u(1-\varepsilon)=0$, where $\varepsilon << 1$.   Table~\ref{table:znum} lists numerical predictions of $z$ for decreasing $\varepsilon$.  For each value of $\varepsilon $, the numerical $\alpha$ step size used in the Runge-Kutta implementation was successively decreased until convergence was established to within the reported digits of the table.  
 
\begin{table}[h!]
\centering
\caption{Predictions of $u(0)\equiv z$ of~(\ref{FP2}), obtained numerically by shooting to the boundary condition $u(1-\varepsilon)=0$. The reported digits are within the convergence tolerance of the shooting iteration. For comparison, the previous benchmark values reported by~\cite{Boyd1997} and~\cite{Maleki} are $z=-2.3919564$ and $-2.391956403$, respectively. }
%\small
\label{table:znum}
\begin{tabular}{cc}
\hline
$\varepsilon$         & $z$    \\ \hline
$10^{-1}$         & -2.6972451986   \\
$10^{-2}$         & -2.4011647700    \\
$10^{-3}$         & -2.3921055641    \\
$10^{-4}$       & -2.3919584270    \\
$10^{-5}$      & -2.3919564287    \\
$10^{-6}$      & -2.3919564035   \\
$10^{-7}$        &-2.3919564032    \\
$10^{-8}$ & -2.3919564032 \\ 
%\cite{Maleki}& -2.391956403  \\
%\cite{Boyd1997} & -2.3919564\\
\hline
\end{tabular}
\end{table}

Using the high precision $z$ value of -2.3919564032 in Table~\ref{table:znum}, the magnitude $S$ of the closest singularity to the origin in the complex $r$-plane is computed by applying the ratio-test to~(\ref{FPseries}) and making a Domb-Sykes plot of the ratios versus inverse coefficient order~\cite{DombSykes}, shown in Fig.~\ref{fig:DS} and taken up to $n=2000$.  A linear fit is then made near $1/(n+1)$=0 with an $n=\infty$ intercept of $S=2.611541077$, where the digits have been conservatively truncated since the residual of the linear fit is of order $10^{-12}$.  
\begin{figure*}[h!]
\centering
\includegraphics[width=3in]{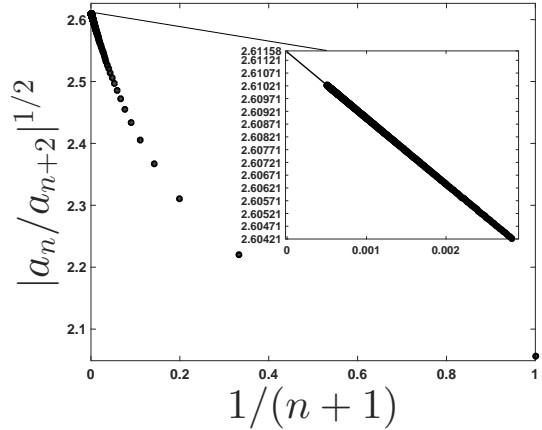}
\caption{Domb-Sykes ratio plot of the coefficients in~(\ref{FPseries}) using $a_0\equiv z=-2.3919564032$ from Table~\ref{table:znum}. The intercept at $1/(n+1)$=0 is 2.611541077.}
\label{fig:DS}
\end{figure*}
\end{document}